\documentclass[11pt]{article}
\usepackage{amssymb,amsmath,enumerate,latexsym, graphicx, epsfig}
\usepackage[latin1]{inputenc}
\usepackage{graphicx}
\usepackage{epsfig}

\def\centerbmp#1#2#3{\vskip#2\relax\centerline{\hbox to#1{\special
  {bmp:#3 x=#1, y=#2}\hfil}}}

\newtheorem{lemma}{Lemma}
\newtheorem{proposition}{Proposition}
\newtheorem{theorem}{Theorem}

\newtheorem{remark}{Remark}

\newtheorem{definition}{Definition}

\newtheorem{assertion}{Assertion}

\newenvironment{proof}{\smallskip\noindent{\it Proof.}\hskip \labelsep}
                        {\hfill\penalty10000\raisebox{-.09em}{$\Box$}\par\medskip}

\def\R{\mathbb{R}}
\def\N{\mathbb{N}}
\def\C{\mathbb{C}}
\def\Z{\mathbb{Z}}
\def\esf{\mathbb{S}}

\def\T{\mathbb{T}}
\def\M{\mathbb{M}}
\def\D{\mathbb{D}}

\def\a{{\alpha}}
\def\be{{\beta}}
\def\t{{\theta}}
\def\g{{\gamma}}
\def\G{{\Gamma}}
\def\l{{\lambda}}

\def\de{{\delta}}

\def\ve{{\varepsilon}}

\def\rth{{\mathbb{R}^3}}
\def\x{{\bf x}}
\def\y{{\bf y}}
\def\Res{{\rm \,Res\,}}

\begin{document}
 \begin{title}
{The classification of doubly periodic minimal tori with parallel ends}
\end{title}
\vskip .5in

\begin{author}
{Joaqu\'\i n P\' erez\thanks{Research partially supported by a MCYT/FEDER grant no.
BFM2001-3318.}, Magdalena Rodr\'\i guez$^*$ \ and Martin Traizet}
\end{author}

\maketitle

\noindent {\sc Abstract.} {\footnotesize Let ${\cal K}$ be the space of properly
embedded minimal tori in quotients of $\R^3$ by two independent translations, with
any fixed (even) number of parallel ends. After an appropriate normalization, we prove that
${\cal K}$ is a 3-dimensional real analytic manifold
that reduces to the finite coverings of the examples defined by Karcher, Meeks and Rosenberg in~\cite{ka4,ka6,mr3}.
The degenerate limits of surfaces in ${\cal K}$ are
the catenoid, the helicoid and three 1-parameter families of surfaces: the simply and doubly periodic Scherk minimal surfaces and the Riemann minimal examples.
}

\section{Introduction}
\label{secintrod} In 1988, Karcher \cite{ka4} defined a $1$-parameter family of
minimal tori in quotients of $\R^3$ by two independent translations.
Each of these surfaces, called {\it toroidal halfplane layer,}  has four
parallel Scherk-type ends in its smallest fundamental domain, is invariant by
reflection symmetries in three orthogonal planes and contains four parallel straight
lines through the ends that project orthogonally onto the corners of a planar
rectangle contained in one of the reflective symmetry planes, and the edges of this rectangle are
just the half-period vectors of the surface, see Figure~\ref{StandardExamples} left. Thanks to this richness of
symmetries, he gave explicitly the Weierstrass representation of the
toroidal halfplane layers in terms of elliptic functions on a $1$-parameter family
of rectangular tori (these examples will be denoted as $M_{\t ,0,0}$ with $0<\t <\frac{\pi }{2}$ in
Section~\ref{secstandardexamples}). Inside a brief remark in his paper and later in
another work~\cite{ka6}, Karcher exposed two distinct $1$-parameter deformations of each
$M_{\t ,0,0}$ by losing some of their symmetries (denoted by
$M_{\t ,\a ,0},M_{\t ,0,\be }$
in Section~\ref{secstandardexamples}).
In 1989, Meeks and
Rosenberg \cite{mr3} developed a general theory for doubly periodic minimal
surfaces with finite topology in the quotient, and used a completely different
approach to find again the examples $M_{\t,0,\be }$ (although up to now it was not clear that Meeks and Rosenberg's examples were the same as Karcher's).
In fact, it is not difficult to produce a 3-parameter family of examples $M_{\t ,\a ,\be }$ containing all the above examples, see Section~\ref{secstandardexamples}.

Hauswirth and Traizet~\cite{HausTraizet1} have extended
previous ideas by P\' erez and Ros~\cite{pro4} and P\' erez~\cite{perez3}
to study the moduli space of all properly embedded doubly periodic minimal
surfaces with a given fixed finite topology in the quotient, proving that in the parallel ends case
(resp. nonparallel) and after identifying by translations, homotheties and rotations,
this moduli space is a real analytic manifold of dimension 3 (resp. 1) around a
{\it nondegenerate} surface.
The nondegeneracy condition concerns the
space of infinitesimal deformations of a given surface; we will see in Section~\ref{secstandardexamples}
that each example $M_{\t ,\a ,\be }$ is nondegenerate.
All these facts point out to a possible global uniqueness of the four ended surfaces $M_{\t ,\a ,\be}$ and their
$k$-sheeted coverings (here $k$ is any positive integer), which in the sequel will be referred to as {\it standard examples,} among all properly embedded doubly periodic minimal surfaces with genus one
and any number of parallel ends in the quotient. In this paper we solve this uniqueness in the affirmative.

\begin{theorem}
\label{thm1}
If $M$ is a properly embedded doubly periodic minimal surface with genus one in the quotient and parallel ends, then $M$ is a standard example.
\end{theorem}

\begin{remark}
The theorem does not hold if we remove the hypothesis on the ends to be parallel, as
demonstrate the 4-ended tori discovered by Hoffman, Karcher and Wei in~\cite{howe3}.
\end{remark}
The analysis of the space of standard examples ${\cal K}$ with $4k$ parallel ends
can be obviously reduced to the case $k=1$ by taking $k$-sheeted coverings. In the four ended case, the
space ${\cal K}_1=\{M_{\t ,\a ,\be}\}_{\t,\a,\be}$ of standard examples is studied in detail in~\cite{PeRo1}.
${\cal K}_1$ is a $3$-dimensional real analytic manifold and the degenerate limits of sequences in ${\cal K}_1$ are the catenoid, the helicoid, any singly or doubly periodic Scherk minimal surface and any Riemann minimal example  (see for instance~\cite{mpr1} for a description
of these last surfaces). Furthermore ${\cal K}_1$
is self-conjugate, in the sense that the conjugate surface\footnote{Two minimal surfaces
$M_1,M_2\subset \R^3$ are {\it conjugate} if the coordinate functions of $M_2$ are harmonic conjugate
to the coordinate functions of $M_1$.} of any element in ${\cal K}_1$ also belongs to ${\cal K}_1$.
Since any standard example admits an orientation reversing involution without
fixed points, Theorem~\ref{thm1} also gives a classification of all properly
embedded minimal Klein bottles with parallel ends in doubly periodic quotients of~$\R^3$.

Theorem~\ref{thm1} can be viewed as a first step toward the
classification of all properly embedded minimal surfaces of genus zero in the
complete flat $3$-manifold $\R^2\times \esf^1$. Concerning this open problem and
except for the flat examples (which reduce to the plane $\R^2\times \{ \theta \} $
and the cylinder $\R \times \esf^1$), all the simple ends of any such a surface are
simultaneously asymptotic to helicoids or to Scherk minimal surfaces.
In the finite topology case, the number of ends is necessarily even and
the helicoid is known to be the unique genus zero example with a
finite number of helicoidal ends (P\' erez and Ros~\cite{pro1}), and it is
conjectured that for any integer $k\geq 2$, the only examples with genus zero and $2k$ ends are the
$(2k-3)$-parameter family of saddle towers having as building blocks the conjugate surface of a Jenkins-Serrin graph
over the right triangle which, after symmetrization, generates a convex $2k$-gon with all edges of the same length (in particular, for $k=2$ the only examples are the singly-periodic Scherk minimal surfaces). The validity of this conjecture will be proved by Pérez and Traizet in a forthcoming paper~\cite{PeTra1}.

The proof of Theorem~\ref{thm1} is a modified application of the machinery developed by Meeks,
P\' erez and Ros in their characterization of the Riemann minimal examples~\cite{mpr1},
hence this reference could be a helpful source to the reader. For $k\in \N$ fixed, one considers the
space ${\cal S}$ of properly embedded doubly periodic minimal surfaces of genus one
in the quotient and $4k$ parallel ends. The goal is to prove that ${\cal S}$
reduces to the space ${\cal K}$ of standard examples (i.e. $k$ sheeted coverings to the surfaces $M_{\t ,\a ,\be }
\in {\cal K}_1$).
The argument is based on modeling ${\cal S}$ as an analytic subset in a complex manifold ${\cal W}$ of finite dimension (roughly, ${\cal W}$ consists of all admissible Weierstrass data for our problem). Then the procedure has three steps:
\begin{itemize}
    \item {\it Properness:} Uniform curvature estimates are proven for a
    sequence of surfaces in ${\cal S}$ constrained to
    certain natural normalizations in terms of the period vector at the ends and
    of the flux of these surfaces (this flux will be defined in Section~\ref{secW}).
    \item  {\it Openness:} Any surface in ${\cal S-K}$ can be minimally deformed by moving its period at
the ends and its flux. This step depends on the properness part and both together imply,
assuming ${\cal S}-{\cal K}\neq \mbox{\O }$ (the proof of Theorem~\ref{thm1} is
by contradiction), that any period at the ends and flux can be achieved by
nonstandard examples.
    \item {\it Uniqueness around a boundary point of ${\cal S}$:}
    Only standard examples can occur nearby a certain minimal surface outside ${\cal S}$ but obtained as a smooth limit
of surfaces in~${\cal S}$. This property together with the last sentence in the
    openness point lead to the desired contradiction, thereby proving Theorem~\ref{thm1}.
\end{itemize}
\noindent
Although the above strategy is adapted from the one in~\cite{mpr1}, there are several major differences between
that work and ours, among which we would like to emphasize two. The main tool in our proof
is a detailed study of the map $C$ that associates
to each $M\in {\cal S}$ two geometric invariants: its period at the ends and
its flux along a nontrivial homology class with vanishing period vector.
The first main difference between this work and~\cite{mpr1} is that $C$ is not
proper, in contrast with the properness of the flux map in~\cite{mpr1}.
Fortunately $C|_{{\cal S}-{\cal K}}$ becomes proper (recall we assumed ${\cal S}-{\cal
K}\neq \mbox{\O }$). This restricted properness follows from curvature estimates as in the first step of the above procedure, together with a local uniqueness argument similar to the third step, performed around any singly periodic Scherk minimal surface, considered as a point in the boundary $\partial {\cal S}$.
In fact, to describe the complete list of limits of sequences in ${\cal S}$ we need a new characterization
of the singly periodic Scherk minimal surfaces
among all properly embedded singly periodic genus zero minimal surfaces with Scherk-type ends, provided that all
ends are parallel except two of them. This result is a special case of the main Theorem in~\cite{PeTra1} and it will be used here.

Once we know that $C|_{{\cal S}-{\cal K}}$ is proper and open, we need a local uniqueness result at a point of
$\partial {\cal S}$ other than a singly periodic Scherk minimal surface, to conclude the third step in our strategy. This
boundary point will be the catenoid, and the local uniqueness follows from
the Inverse Function Theorem. The second main discrepancy with~\cite{mpr1} is a technical
difficulty: nearby the catenoidal limit, the ends of surfaces in ${\cal S}$ group in couples, each
couple giving rise to a single end of the limit catenoid. Such a collapsing
phenomenon makes certain residues to blow-up, and a careful study of speeds of
degeneration is needed to rectify the mapping to which the Inverse Function Theorem
is applied to.

The paper is organized as follows. In Section~\ref{secprelim} we recall the
necessary background to tackle our problem. Section~\ref{secstandardexamples} is
devoted to introduce briefly the 3-parameter family ${\cal K}$ of standard examples.
The complex manifold of admissible
Weierstrass data ${\cal W}$ and the natural mappings on it are
studied in Section~\ref{secW}. In Section~\ref{seccurvestim} we obtain the
curvature estimates needed for the first point of our strategy.
Sections~\ref{secuniqaroundScherk1p} and~\ref{secuniqaroundcat} deal with the local
uniqueness around the singly periodic Scherk minimal surfaces and the catenoid,
respectively. The second point of our above strategy (openness) is the goal of
Section~\ref{secopen}, and finally Section~\ref{secproofmainthm} contains the proof of Theorem~\ref{thm1}.

The authors would like to thank Antonio Ros and Bill Meeks by helpful conversations.

\section{Preliminaries.}
\label{secprelim}
Let $\widetilde{M}\subset \R^3$ be a connected orientable\footnote{Unless otherwise stated, all surfaces in the paper are supposed to be orientable.}
properly embedded minimal surface,
invariant by a rank 2 lattice ${\cal P}$ generated by two linearly independent translations
$T_1, T_2$ (we will shorten by calling $\widetilde{M}$ a {\it doubly periodic minimal
surface}). $\widetilde{M}$ induces a properly embedded minimal surface $M=\widetilde{M}/{\cal
P}$ in the complete flat $3$-manifold $\R^3/{\cal P}= \T^2\times \R $, where $\T $
is a $2$-dimensional torus.
Reciprocally, if $M\subset \T\times \R$ is a properly embedded nonflat minimal
surface, then its lift $\widetilde{M}\subset\R^3$ is a connected doubly
periodic minimal surface by the Strong Halfspace Theorem of Hoffman and
Meeks~\cite{hm10}. Existence and classification theorems for doubly periodic
minimal surfaces are usually tackled by considering the quotient surfaces in $\T
\times \R $. An important result by Meeks and Rosenberg~\cite{mr3} insures that
a properly embedded minimal surface $M \subset \T\times \R$ has finite topology if
and only if it has finite total curvature, and in this case $M$ has an even number
of ends, each one asymptotic to a flat annulus (Scherk-type end). Later, Meeks~\cite{me21} proved that
any properly embedded minimal surface in $\T \times \R$ has a finite number of ends, so the finiteness of its
total curvature is equivalent to the finiteness of its genus.

When normalized so that the lattice of periods ${\cal P}$ is horizontal, we distinguish two types of ends, depending
on whether the well defined third coordinate function on $M$ tends to $\infty$ ({\it top} end) or to $-\infty$ ({\it bottom} end) at the corresponding puncture.
By separation properties, there are an even number of top (resp. bottom) ends.
Because of embeddedness, top (resp. bottom) ends are always parallel each other.
If the top ends are not parallel to the bottom ends, then there exists an algebraic obstruction on the period lattice,
which must be {\it commensurable} as in the classical doubly periodic minimal
surfaces defined by Scherk in~\cite{sche1} or in the 4-ended tori found by Hoffman, Karcher and Wei~\cite{howe3}.
If the top ends are parallel to the bottom ends, then the cardinals of both families of ends coincide, therefore the
total number of ends of $M$ is a multiple of four. For details, see~\cite{mr3}.

We will focus on the parallel ends setting, where the simplest possible topology is a finitely punctured torus
(properly embedded minimal planar domains in $\T \times \R $ must have nonparallel ends by Theorem 4 in~\cite{mr3}; in fact Lazard-Holly and Meeks~\cite{lhm} proved that the doubly periodic Scherk minimal surfaces are the unique possible examples with genus zero). Our main goal is to give a complete classification of all examples with genus one and parallel ends. To do this, we will normalize appropriately the surfaces under consideration.

Given a positive integer $k$, let
${\cal S}$ be the space of all properly embedded minimal tori in a quotient $\R^3/{\cal P}=\T\times\R$ modulo a rank 2 lattice ${\cal P}$ generated by two independent translations (which depend on the surface), one of them being in the direction of the $x_2$-axis, with $4k$ horizontal Scherk-type ends.
Given $M\in\cal S$, we denote respectively by $P_{\G}\in {\cal P}$ and $F_\G$ the period and flux vectors of $M$ along an oriented closed curve $\G \subset M$.
By the Divergence Theorem, $P_{\G }$, $F_\G$ only depend on the homology class of $\G $ in $M$.
The period and flux vectors
$H$, $ F$ at an end of $M$ (defined as the period and flux along a small loop
around the puncture with the inward pointing conormal vector respect to the disk that
contains the end) are related by the equation $ F=H\wedge  N_0$, where $N_0$ is the
value of the $\esf^2$-Gauss map at the puncture. In our normalization, each of the period vectors at the ends of $M$ is of the form $H=\pm (0,\pi a,0)$ with $a>0$.
The end is called a {\it left end} if $F=(-\pi a,0,0)$, and a {\it right end} if $F=(\pi a,0,0)$.
As $M$ is embedded, each family of ``sided'' ends is naturally ordered by heights;
in fact the maximum principle at infinity~\cite{mr1}
implies that consecutive left (resp. right) ends are at positive
distance. Furthermore, their limit normal vectors are opposite by a trivial separation argument.

We will denote by $\widetilde M\subset \R^3$ the (connected) doubly periodic
minimal surface obtained by lifting $M$. Since points in $\widetilde{M}$ homologous by ${\cal P}$
have the same normal vector, the stereographically projected Gauss map $g:\widetilde{M}\to
\overline{\C}=\C \cup \{\infty\} $
descends to $M$. As $M$ has finite total curvature, $g$ extends meromorphically to
the conformal torus $\M$ obtained after attaching the ends to $M$, with values $0,\infty $
at the punctures. As the period lattice ${\cal P}$ is not horizontal,
the third coordinate function $x_3$ of $\widetilde{M}$ is multivalued on $M$ but the
{\it height differential} $dh =\frac{\partial x_3}{\partial z}dz$ defines a univalent
meromorphic differential on $M$ (here $z$ is a holomorphic coordinate). Since $M$ has finite total curvature and horizontal ends,  $dh$ extends to
a holomorphic differential on $\M $. The next statement collects some elementary properties
of any surface in ${\cal S}$. Given $v\in {\cal P}-\{0\} $,
$\widetilde{M}/v$ will stand for the singly periodic minimal surface obtained as the
quotient of $\widetilde{M}$ by the translation of vector $v$.

\begin{proposition}
\label{propos1}
Given $M\in {\cal S}$, it holds
\begin{description}
\item[{\it 1.}]$g:\M \to \overline{\C }$ has degree $2k$, total branching number $4k$, does not take vertical directions on $M$ and it is unbranched at the ends.
\item[{\it 2.}] The period vectors at the ends coincide up to sign and we will denote them by $H=\pm(0,\pi a,0)$, $a>0$.
The period lattice ${\cal P}$ of $\widetilde{M}$ is generated by $H$ and a nonhorizontal vector $ T\in \R^3$, this last one being the period vector along a closed curve $\g _1\subset M$
such that $[\g _1]\neq 0$ in the homology group $H_1(\M ,\Z)$.
\item[{\it 3.}] Let ${\cal E}$ be the set of Scherk-type ends of $\widetilde{M}/H$. Then
$(\widetilde{M}/ H)\cup{\cal E}$ is conformally $\C^*=\C-\{0\}$, and the height differential writes as $dh=c \frac{dz}{z}$ in $\C^*$, with $c\in \R^*=\R -\{0\}$.
\end{description}
Let $\Pi \subset \R^3$ be a horizontal plane.
\begin{description}
\item[{\it 4.}] If $\Pi /H$ is not asymptotic to an end in ${\cal E}$, then
$(\widetilde M\cap \Pi )/ H$ is transversal and connected. The period vector along $(\widetilde M\cap \Pi )/ H$ either
vanishes or equals $\pm H$.
\item[{\it 5.}] We divide ${\cal E}$ in {\it right ends} and {\it left
ends,} depending on whether the flux vector at the corresponding end (with the inward pointing conormal vector) is $(a,0,0)$ or $(-a,0,0)$, respectively. If $\Pi /H$ is asymptotic to an end in ${\cal E}$, then
$(\widetilde M\cap \Pi )/ H$ consists of one or two properly
embedded arcs. If $(\widetilde M\cap \Pi )/ H$ reduces to one arc $\G $, then both ends of $\G $
diverge to the same end in ${\cal E}$. In the two arcs intersection case, both arcs
travel from one left end to one right end in ${\cal E}$.
\item[{\it 6.}] There exists an embedded closed curve $\g_2\subset M$ such that $\{[\g_1],[\g_2]\}$ is basis of $H_1(\M , \Z )$ and $ P_{\g _2}=0$. Up to orientation, $\g _2$ represents the
 unique nontrivial homology class in $H_1(\M,\Z )$ with associated period zero and an embedded
 representative.
\item[{\it 7.}] Let $[\g ]\in H_1(M,\Z )$ be a homology class with an embedded representative that generates the
homology group of $(\widetilde{M}/H)\cup{\cal E}$. Then
the third component $(F_\g)_3$ of the flux of $M$ along any representative $\g\in [\g ]$ neither vanish nor depends
on $[\g ]$ (up to orientation).
\end{description}
\end{proposition}
\begin{proof}
Since $dh$ has no poles on the torus $\M $, it cannot have zeros, which implies
that the only zeros and poles of $g$ are at the ends. As at each puncture one of the meromorphic differentials $g\, dh,\frac{dh}{g}$ has a simple pole (Meeks and Rosenberg~\cite{mr3}), we conclude that $g$ is unbranched at the ends.
Since $M$ has $4k$ ends, $g$ must have degree $2k$ and the Riemann-Hurwitz formula implies that its total branching number is $4k$, so statement {\it 1} is proved.

Denote by $T_1, T_2$ two generators of $\cal P$. As $M$ is properly embedded with
horizontal ends,  we can assume that $T_1$ is not horizontal and $T_2$ is horizontal.
It follows that the period vector at an end
$E$ of $M$ is $H=n T_2$ with $n\in\Z-\{ 0\} $, and that $T_1$ is the period vector
along a closed curve $\g_1\subset M$ with $[\g _1]\neq 0$ in $H_1(\M ,\Z)$.
Part~{\it 2} of the Proposition will be proved if we show that $n=\pm 1$, by taking $T=T_1$.
Since $E$ is invariant by $T_2$ and the period vector of $E$ is $nT_2$, it follows that $g$ has
branching order  $|n|-1$ at the puncture. As $g$ is unbranched at the ends, we deduce that $n=\pm 1$
as desired.

Since $M$ is a finitely punctured torus, the singly periodic surface $\widetilde{M}/H$ is a
cylinder with infinitely many points removed. Note that the third coordinate function $x_3$
is a well defined harmonic function on $\widetilde{M}/H$ that extends smoothly through ${\cal E}$, giving rise to a proper harmonic function
on $(\widetilde{M}/H)\cup {\cal E}$. Therefore, $(\widetilde{M}/H)\cup {\cal E}$
is conformally $\C^*=\C-\{0\}$ and $x_3$ writes as $x_3(z)=c\ln |z|+c'$  with $c\in \R^*$, $c'\in \R $,
which is statement~{\it 3} of the Proposition.
This description of $x_3$ implies that $(\widetilde M\cap \Pi )/H$
corresponds in the $\C^*$-model to a
possible punctured circle $C_r=\{|z|=r\} $ for certain $r>0$. The
hypotheses in item~{\it 4} correspond to the case that $C_r$ does not
contain ends in ${\cal E}$, so the conclusion of {\it 4} is clear. Under the
hypotheses of item~{\it 5}, $C_r$ contains at most two ends of ${\cal E}$ (one left end
and/or one right end because ends of the same side are separated by heights), and
{\it 5} also holds easily.

In order to see item {\it 6}, let $\be $ be a compact horizontal level section of $M$.
By the description above, $\be $ generates the homology of $(\widetilde{M}/H)\cup {\cal E}=\C ^*$.
Since the period vector $T=P_{\g _1}$ is not horizontal, we conclude that
$\{ [\g_1],[\be ]\} $ is a basis of $H_1(\M,\Z )$. If $ P_{\be }=0$, then the first assertion in {\it 6}
is proved with $\g _2=\be $. If $P_{\be }\neq 0$, then $ P_{\be }=\pm  H$ by item~{\it 4}. In this case
we choose as $\g _2$ an embedded closed curve in $M$ homologous to $\be $ in $\C^*$ such that $\be \cup \g _2$
bounds just an end with period vector $-P_{\be }$. Finally, suppose that $[\G ]\in H_1(\M,\Z )-\{ 0\} $ has
an embedded representative $\G $ (which can be assumed to lie in $M$) with $ P_{\G}=0$.
Since $\g _1,\g _2$ and small loops $\a _1,\ldots ,\a _{4k}$ around the punctures generate $H_1(M,\Z )$, we can write
\begin{equation}
\label{eq:lema1A} [\G ]=a_1[\g _1]+a_2[\g _2]+\sum _{i=1}^{4k}b_i[\a _i] \quad \mbox{in
}H_1(M,\Z )
\end{equation}
for integers $a_1,a_2,b_1,\ldots ,b_{4k}$. Taking periods in (\ref{eq:lema1A}) and
having in mind that $ P_{\g _1}= T$ and that the periods at the ends are $\pm  H$, we
obtain $0= P_{\G}=a_1  T+b H$, where $b\in \Z $. As $ T, H$ are linearly independent, it
follows $a_1=0$. Now (\ref{eq:lema1A}) implies that $[\G ]=a_2[\g _2]$ in $H_1(\M ,\Z
)$, so the embeddedness of $\G $ forces $a_2$ to be $\pm 1$. This proves {\it 6}.

Finally, recall we have shown that any compact horizontal level section $\be $ of $M$ is an embedded closed curve
such that $\{[\g_1],[\be ]\}$ is a basis for $H_1(\M ,\Z )$. As the conormal vector to $M$ along $\be $ has third coordinate with constant sign, we have $(F_\be )_3\neq 0$. Note that if $\g \subset M$ is an  embedded closed curve homologous to $\be $ in $(\widetilde{M}/H)\cup {\cal E}$, then $\g \cup \be $ bounds
a finite number of ends in ${\cal E}$, whose fluxes are all horizontal. Thus $(F_\g )_3=(F_\be)_3$ and the Proposition is proved.
\end{proof}

\begin{remark}
\label{remark1}
\begin{description}
\item[]
\item[{\it (i)}] In general, we cannot expect the curve $\g_2$ in item {\it 6} above to be a compact horizontal section.
For instance, all the horizontal level curves of the standard example $M_{\t ,\frac{\pi}{2},0}$ in Section~\ref{secstandardexamples} are open arcs when viewed in $\R^3$, see Figure~\ref{StandardExamples2} right.
\item[{\it (ii)}] With the notation in Proposition~\ref{propos1}, the fact that all the fluxes at the ends of
$M\in {\cal S}$ point to the $x_1$-axis implies that the component $(F_\g )_2$ in the direction
of the $x_2$-axis is also independent (up to orientation) of the homology class $[\g ]\in H_1(M,\Z )$ satisfying the hypotheses of item~{\it 7}.
\end{description}
\end{remark}
\par
\noindent
There is a last natural normalization on the surfaces in ${\cal S}$, which we now explain. Given $M\in {\cal S}$,
Proposition~\ref{propos1} gives a nontrivial homology class in $H_1(\M ,\Z )$ with an embedded representative
$\g_2\subset M$ such that $P_{\g_2}=0$ and $(F_{\g_2})_3>0$. In the sequel, we will always normalize our surfaces so that  $(F_{\g_2})_3=2\pi $, which can be achieved after an homothety. Note that this normalization
is independent of the homology class of $\g _2$ in $H_1(M ,\Z )$ (up to orientation), see item~{\it 7} of Proposition~\ref{propos1}.

We label by $\widetilde{\cal S}$ the set of {\it marked surfaces} $(M,p_1,\ldots ,p_{2k},q_1,\ldots ,q_{2k},[\g _2])$
where
\begin{enumerate}
\item $M$ is a surface in ${\cal S}$ whose period lattice is generated by $H,T\in \R^3$, where $H=(0,a,0)$, $T=(T_1,T_2,T_3)$ and $a,T_3>0$;
\item $\{ p_1,\ldots ,p_{2k}\} =g^{-1}(0)$, $\{ q_1,\ldots ,q_{2k}\} =g^{-1}(
\infty )$ and the ordered lists $(p_{1},q_{1},\ldots ,p_{k},q_{k})$, $(p_{k+1},q_{k+1},\ldots ,p_{2k},q_{2k})$ are the two families of ``sided'' ends of $M$, both ordered by increasing heights in the quotient;
\item $[\g _2]\in H_1(M,\Z )$ is the homology class of an embedded closed curve $\g _2\subset M$ satisfying $P_{\g _2}=0$, $(F_{\g_2})_3=2\pi $. We additionally impose that $\g_2$ lifts to a curve contained in a fundamental domain of the doubly periodic lifting of $M$ lying between two horizontal planes $\Pi ,\Pi +T$.
\end{enumerate}
We will identify in $\widetilde{\cal S}$ two marked surfaces that differ by a translation that preserves both orientation, the above ``sided'' ordering of their lists of ends and the associated homology classes. The same geometric surface
in ${\cal S}$ can be viewed as different marked surfaces in $\widetilde{\cal S}$.
We will simply denote as $M\in \widetilde{\cal S}$ the marked surfaces unless it leads to confusion.
$\widetilde{\cal S}$ can be naturally endowed with a topology: a sequence of marked surfaces $\{ M_n\} _n\subset \widetilde{\cal S}$ converges to $M\in \widetilde{\cal S}$ if the associated sequence of minimal surfaces $\{ M_n\} _n\subset {\cal S}$ converges smoothly to $M\in {\cal S}$ (in the uniform topology on compact sets),
the ordered list of ends associated to $M_n$ converges to the corresponding one for $M$ and the homology classes
$[\g _{2,n}]\in H_1(M_n,\Z )$ in the last component of the marked surfaces $M_n$ have representatives converging uniformly to a representative of the last component of $M$. With this topology, a geometric minimal surface $M\in {\cal S}$ produces a finite subset in $\widetilde{\cal S}$.

Consider $M\in \widetilde{\cal S}$ with Gauss map $g$ and height differential $dh$.
An elementary calculation gives the periods $P_{p_j},P_{q_j}$ and fluxes $F_{p_j},F_{q_j}$ at the ends of $M$ as follows:
\begin{equation}
\label{eq:periodfluxends}
P_{p_j}+i\, F_{p_j}=\pi \mbox{Res}_{p_j}\left( g^{-1}\, dh \right) (i,-1,0),\qquad
P_{q_j}+i\, F_{q_j}=-\pi\mbox{Res}_{q_j}(g\, dh)(i,1,0),
\end{equation}
where Res$_A$ denotes the residue of the corresponding meromorphic differential at a point~$A$.
The fact that $P_{p_j},P_{q_j}$ point to the $x_2$-axis translates into $\mbox{Res}_{p_j}
\left( g^{-1}\, dh\right) $, $\mbox{Res}_{q_j}(g\, dh)\in \R $. By definition of the ordering of the ends of $M$ as a
marked surface, we have that
\begin{equation}
\label{eq:gcloseends}
\mbox{Res}_{p_j}\left( g^{-1}\, dh \right) =
-\mbox{Res}_{q_j}( g\, dh)=\left\{
\begin{array}{rl}
a & (1\leq j\leq k)\\
-a & (k+1\leq j\leq 2k),
\end{array}\right. \qquad
\end{equation}
for certain $a\in\R^*$ (the case $a>0$ corresponds to $p_1,\ldots ,p_k,q_1,\ldots ,q_k$ being
right ends of~$M$). Recall that $P_{\g _2}=0$ and $(F_{\g_2})_3=2\pi $. Thus,
\begin{equation}
\label{eq:gcloseshomology}
\overline{\int_{\g_2}g^{-1}\, dh}=\int_{\g_2}g\, dh ,\qquad \int _{\g _2}dh=2\pi i.
\end{equation}

\section{Standard examples.}
\label{secstandardexamples}

We dedicate this Section to introduce briefly the $3$-parameter family of standard examples ${\cal K\subset S}$ to which the uniqueness Theorem~\ref{thm1} applies. Some of the properties of ${\cal K} $ in this Section are
long but straightforward computations that can be found in detail in~\cite{PeRo1}.
As the standard examples with $4k$-ends are nothing but $k$-sheeted coverings of $4$-ended standard examples in ${\cal K}_1$, we will concentrate on these last ones. Each $M_{\t,\a,\be}\in{\cal K}_1$ is determined by the $4$ branch values of its Gauss map, which consist of two pairs of antipodal points $D,D',D''=-D,D'''=-D'$ in the sphere $\esf^2$.
Since the Gauss map of any surface in ${\cal S}$ is unbranched at the ends, $D,D',D'',D'''$ must be different from the North and South Poles. We also let $e$ be the equator in $\esf^2$ that contains
$D,D',D'',D'''$. Given a point $P\in e$, the branch values of the Gauss map can be determined by giving only one angle
$\t \in (0,\frac{\pi }{2})$, in such a way that the position vectors of $D,D'$ form an angle of $2\t $ and the position vector of $P$ bisects such an angle. We will call a {\it spherical configuration} to any set $\{D,D',D'',D'''\}$ as above.

\begin{figure}
\centerline{\includegraphics[width=13.6cm,height=4.64cm]{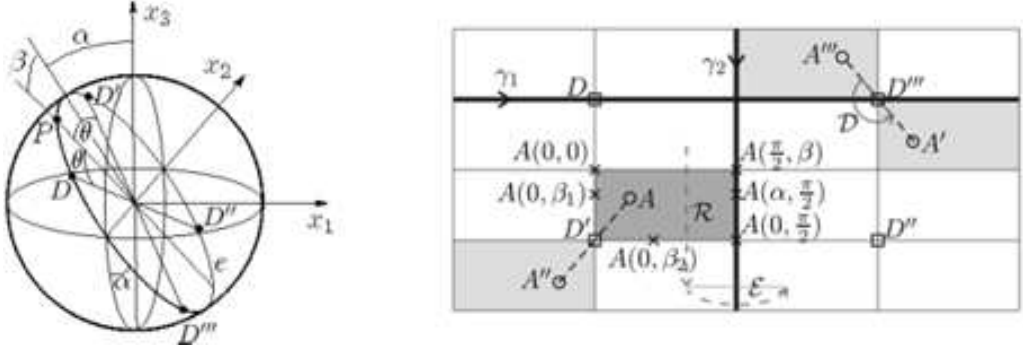}}
\caption{Left: Spherical configuration of $M_{\t ,\a ,\be }$. Right: Behavior of the
pole $A(\a,\be )$ of $g$ in the dark shaded rectangle ${\cal R}$ (here,
$0<\be_1<\t<\be_2<\frac{\pi}{2}$, $\a\in(0,\frac{\pi}{2})$ and
$\be\in[0,\frac{\pi}{2}]$). The remaining ends of $M_{\t ,\a ,\be }$ move in the
light shaded rectangles as $\a ,\be $ vary. } \label{figure8}
\end{figure}

Given $(\t ,\a ,\be) \in (0,\frac{\pi }{2})\times [0,\frac{\pi }{2}]\times [0,\frac{\pi }{2}]$ with
$(\a ,\be )\neq (0,\t )$, we define the spherical configuration of the potential standard example $M_{\t ,\a ,\be }$
by means of the angle $\t $, the equator $e$ and the point $P$ as follows:
\begin{enumerate}
\item Let $e_0$ be the inverse image of the imaginary axis in $\C $ through the stereographic projection from the North Pole of $\esf ^2$. Then $e$ is the image of $e_0$ through the rotation by angle $\a $ around the $x_2$-axis.
\item $P$ is the image of the North Pole through the composition of a rotation by angle $\be $ around the $x_1$-axis
with a rotation by angle $\a $ around the $x_2$-axis, see Figure~\ref{figure8} left.
\end{enumerate}

\begin{figure}
\centerline{\includegraphics[width=13.6cm,height=4cm]{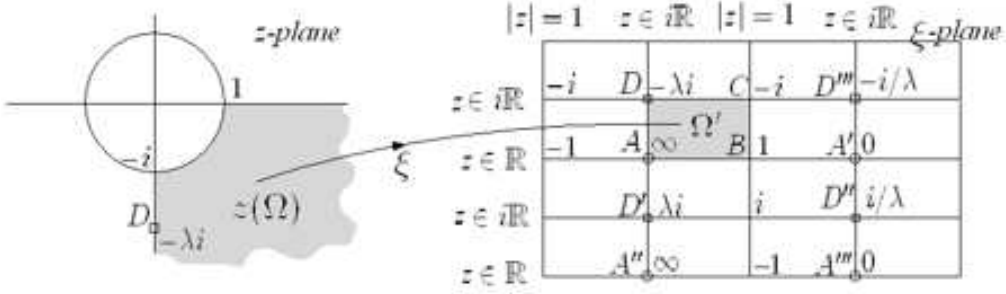}}
\caption{The
biholomorphism $\xi$ between the shaded regions. We have included some values of the
$z$-projection in the $\xi$-plane model of $\Sigma_\t$.}
\label{figA}
\end{figure}

Projecting stereographically from the North Pole the branch points $D,D',D'',D'''$
for the spherical configuration associated to $(\t ,0,0)$, one just finds the four
roots of the polynomial $(z^2+\l^2)(z^2+\l^{-2})$ where
$\l=\l(\t)=\cot\frac{\t}{2}$. Therefore, the underlying conformal compactification
of the potential surface $M_{\t ,0,0}$ is the rectangular torus
\[
\Sigma_\t=\big\{ (z,w)\in \overline{\C }^2\ | \ w^2=(z^2+\l^2)(z^2+\l^{-2})\big\} ,
\]
and its extended Gauss map is the $z$-projection $(z,w)\in \Sigma _\t\mapsto z\in \overline{\C }$ on $\Sigma _\t$.
Note that the spherical configuration for angles $(\t ,\a ,\be )$ differs from the one associated to $(\t ,0,0)$
in a Möbius transformation $\varphi $. Thus the compactification of any $M_{\t ,\a ,\be }$, which is the branched covering of $\esf^2$ through its Gauss map, is $\Sigma _{\t }$.
Furthermore, the composition of the Gauss map of $M_{\t ,0,0}$ with $\varphi $ gives the Gauss map of the potential example $M_{\t ,\a ,\be }$:
\[
g(z,w)=\frac{z\left( i\cos (\frac{\a -\be }{2})+\cos (\frac{\a +\be }{2})\right) +\sin (\frac{\a -\be }{2})+i
\sin (\frac{\a +\be }{2})}{\cos (\frac{\a -\be }{2})+i\cos (\frac{\a +\be }{2})-z\left( i\sin (\frac{\a -\be }{2})
+\sin (\frac{\a +\be }{2})\right) },\quad (z,w)\in \Sigma_\t ,
\]
and its ends are $\{ A,A',A'',A'''\} =g^{-1}(\{ 0,\infty \})$. As the height differential $dh$ of $M_{\t ,\a ,\be }$
is a holomorphic $1$-form on $\Sigma _{\t }$, we have $dh=\mu \frac{dz}{w}$ for certain $\mu =\mu (\t ,\a ,\be )\in \C^*$. It will be also useful to have a second representation of $\Sigma_\t $ as a quotient of the $\xi $-plane $\C$ by two orthogonal translations. Let $\Omega \subset \Sigma_\t $ be one of the two connected components of $\{ (z,w)
\in \Sigma _\t \ / \ |z|>1, -\frac{\pi }{2}<\arg (z)<0\} )$. $\Omega$ is topologically a disk and its boundary contains just one branch point $D_1=(-\l i,w(-\l i))$ and one pole $A_1=(\infty ,\infty )$ of the $z$-projection. Let $\Omega'$  be an open rectangle of consecutive
vertices $A,B,C,D\in \C$ with the segment $\overline{AB}$ being horizontal, such that there exists a biholomorphism
$\xi :\{ |z|>1,-\frac{\pi }{2}<\arg (z)<0\} \to \Omega'$ with boundary values $\xi (\infty )=A$, $\xi (1)=B$,
$\xi (-i)=C$ and $\xi (-\l i)=D$. Then the composition of the $z$-projection with $\xi$ defines a biholomorphism between $\Omega$ and $\Omega'$. After symmetric extension of this biholomorphism across the boundary
curves of $\Omega ,\Omega '$ we will get a biholomorphism from $\Sigma_\t $ to the quotient of the $\xi $-plane modulo
the translations given by four times the sides of the rectangle $\Omega'$, see Figure~\ref{figA}.
Note that if $\a =\be =0$, then $g(z,w)=z$ and the equalities $\xi (\infty )=A$, $\xi (-\l i)=D$ justify to use the same symbols $A,D$ previously defined for two of the corners of the rectangle $\Omega '$.

The above identification allows us to see the ends $A,A',A'',A'''$ as functions of $(\a ,\be )\in[0,\frac{\pi }{2}]^2-\{ (0,\t )\} $ valued in the $\xi $-plane model of $\Sigma_\t$.
$A(\a,\be)$ moves on the dark shaded rectangle ${\cal R}$ in Figure~\ref{figure8} right.
For $(\a ,\be )=(0,0)$, we choose $A(0,0)$ as the up-left corner of~${\cal R}$.
As $\a $ increases in $[0,\frac{\pi }{2}]$, the end $A(\a ,0)$ of $M_{\t ,\a ,0}$ moves continuously and horizontally to the right, until reaching for $\a =\frac{\pi }{2}$ the
up-right corner of ${\cal R}$.
The curve $\be \in [0,\t )\mapsto A(0,\be )$ parametrizes downwards the left vertical boundary edge of ${\cal R}$, until reaching the branch point $D'$.
The curve $\be \in (\t ,\frac{\pi }{2}]\mapsto A(0,\be )$ parametrizes the lower horizontal edge of $\partial {\cal R}$ from $D'$ until the down-right corner of ${\cal R}$.
The curve $\a \in [0,\frac{\pi }{2}]\mapsto A(\a ,\frac{\pi }{2})$ parametrizes upwards the right vertical edge of $\partial {\cal R}$. Finally, the curve $\be \in [0,\frac{\pi }{2}]\mapsto A(\frac{\pi }{2},\be )$ is constant at the up-right corner of ${\cal R}$. These boundary values of $A(\a,\be)$ extend continuously and bijectively from $(0,\frac{\pi }{2})^2$  to the interior of ${\cal R}$.
The behavior of the remaining three ends of $M_{\t ,\a ,\be }$ on the $\xi $-plane model can be deduced from
$A(\a ,\be )$ by using the isometry group Iso$(\t ,\a ,\be )$ of the induced metric $ds^2=\frac{1}{4}(|g|+|g|^{-1})^2|dh|^2$, which we now investigate.

First note that the identity in $\esf^2$ lifts via $g$ to two different isometries of $ds^2$, namely the identity in $\Sigma _{\t }$ and the deck transformation ${\cal D}(z,w)=(z,-w)$, both restricted to $\Sigma _{\t }-g^{-1}(\{ 0,
\infty \} )$. ${\cal D}$ corresponds in the $\xi $-plane to the $180^\circ$-rotation about any of the branch points of
the $z$-projection. The antipodal map $\aleph $ on $\esf^2$ also leaves invariant both the spherical configuration of
$M_{\t ,\a ,\be }$ and the set $\{ (0,0,\pm 1)\} $, and the equality $\left[ w\left( \frac{-1}{\overline{z}}\right) \right] ^2=\frac{[w(z)]^2}{\overline{z}^4}$ for any $(z,w)\in \Sigma _{\t }$ implies that $\aleph $ lifts through $g$ to two isometries of $ds^2$, which we call ${\cal E}$ and ${\cal F}={\cal D}\circ {\cal E}$. Both ${\cal E},{\cal F}$ are antiholomorphic involutions of $\Sigma _{\t }$ without fixed points. This property implies that one of these
involutions, say ${\cal E}$, corresponds in the $\xi $-plane model to the composition of the translation by the vector $\overline{D,D'}$ with the symmetry with respect to the right vertical edge of $\partial {\cal R}$, see Figure~\ref{figure8} right.
The remaining ends of $M_{\t ,\be ,\a }$ in terms of $A=A(\a,\be)$ are (up to relabeling)
\begin{equation}
\label{eq:A',A'',A'''}
A''={\cal D}(A),\qquad A'''={\cal E}(A),\qquad A'={\cal D}(A''').
\end{equation}

We now study the period problem for $M_{\t ,\a ,\be }$. Recall that $dh=\mu
\frac{dz}{w}$ on $\Sigma _{\t }$ with $\mu \in \C^*$. From now on, take $\mu \in \R
^*$. The period $P_A$ and flux $F_A$ of $M_{\t ,\a ,\be }$ at the end $A$ with
$g(A)=\infty $ are given by
\begin{equation}
\label{eq:periodendsnosim}
P_{A}=\pi \mu \big( i\, E(\t,\a,\be),0\big) ,\qquad
F_A=\pi \mu \big( E(\t,\a,\be),0\big) ,
\end{equation}
where we have used the identification $\R^3\equiv \C \times \R$ by $(a,b,c)\equiv (a+ib,c)$, and
 $E(\t,\a,\be)=[\cos^2\a +\csc^2\t (\sin \a \cos \be -i\sin \be )^2]^{-1/2}$ (we have
chosen a branch of $w$ for computing (\ref{eq:periodendsnosim}), which only affects the result up to sign).
The periods and fluxes at $A',A'',A'''$ can be easily obtained using (\ref{eq:A',A'',A'''}), (\ref{eq:periodendsnosim})
and the pullbacks by ${\cal D},{\cal E}$ of the Weierstrass form $\Phi =\big(\frac{1}{2}(g^{-1}-g),\frac{i}{2}(g^{-1}+g),1\big) dh$
of $M_{\t ,\a ,\be }$:
\begin{equation}
\label{eq:D,Epullbacks}
{\cal D}^*\Phi =-\Phi ,\qquad {\cal E}^*\Phi =-\overline{\Phi }.
\end{equation}
Concerning the period problem in homology, let $\g_1,\g_2$ be the simple closed curves in $\Sigma_\t $ obtained
respectively as quotients of the horizontal and vertical lines in the $\xi $-plane passing through $D,D'''$ and through the right vertical edge of $\partial{\cal R}$ (see Figure~\ref{figure8} right). Clearly $\{[\g_1],[\g_2]\}$ is a basis of $H_1(\Sigma_{\t},\Z )$.
Both $\g _1,\g _2$ miss the ends $A,A',A'',A'''$ except for certain extreme values of $(\a ,\be )$ (for instance,
$A(\pi /2,0)$ lies in $\g _2$), but one can think of the homology classes $[\g _1],[\g _2]
\in H_1(M_{\t ,\a ,\be },\Z )$ as being independent of $(\a ,\be )$. Furthermore,
\begin{equation}
{\cal E}_*[\g_1]=-[\g_1]-[\g_A]-[\g_{A'}],\qquad {\cal E}_*[\g_2]=[\g_2] \quad \mbox{ in }H_1(M_{\t ,\a ,\be },\Z ),
\label{eq:calEenhomologia}
\end{equation}
where $\g_A,\g _{A'}$ denote small loops around $A,A'$. Equations (\ref{eq:A',A'',A'''}), (\ref{eq:D,Epullbacks}) and the first equality in (\ref{eq:calEenhomologia}) imply
$-\overline{\int_{\g_1}\Phi }=-\int_{\g_1}\Phi -\int_{\g_A}\Phi +\overline{\int_{\g_A}\Phi }$.
Taking imaginary parts, we find
\begin{equation}
\label{eq:fluxg1Mtab}
F_{\g_1}=-F_{\g_A}=-F_A.
\end{equation}
Similarly, the second equalities in (\ref{eq:D,Epullbacks}), (\ref{eq:calEenhomologia}) insure that $\int_{\g_2}\Phi =
-\overline{\int_{\g_2}\Phi }$, whose real part gives
\begin{equation}
\label{eq:Perg2Mtab}
P_{\g_2}=0.
\end{equation}
Equation (\ref{eq:Perg2Mtab}) implies that we can take $\g _2$ as the embedded closed curve appearing in
item~{\it 6} of Proposition~\ref{propos1} (except for the extreme values of $(\a ,\be )$ mentioned above, in which we deform $\g _2$ keeping $[\g _2]\in H_1(M_{\t ,\a ,\be },\Z )$ constant in $\a,\be$). As $dh$ is holomorphic and nontrivial on
$\Sigma _{\t }$, (\ref{eq:Perg2Mtab}) also implies that $\int_{\g_2}dh\in i\R^*$, so we must rescale to have
$\int _{\g _2}dh=2\pi i$ as part of the normalization of our surfaces in ${\cal S}$. This rescaling gives $\mu $:
\begin{equation}
\label{eq:mu}
\mu=\mu (\t )=\frac{\pi }{2{\displaystyle \int_{-\pi /2}^0\sqrt{2\cos (2t)+\l^2(\t )+\l^{-2}(\t )}\, dt}}
=\frac{\pi }{\sin \t \, \mbox{EllipticK}(\sin^2\t )},
\end{equation}
where EllipticK$(m)=\int_0^{\pi /2}\sqrt{1-m\sin^2u}\, du$ is the complete elliptic integral of the first kind.
Since $\Re \int _{\g _2}dh=0$ and  $dh$ is holomorphic on $\Sigma_\t $, we deduce that $P_{\g_1}$ has nonvanishing third component.
In particular, $P_{A}$ and $P_{\g_1}$ are linearly independent.
All these facts imply that $M_{\t ,\a ,\be }$ is a complete immersed doubly periodic minimal torus with four horizontal embedded Scherk-type ends and period lattice generated by $P_{A},P_{\g_1}$.

We claim that $M_{\t ,\a ,\be }$ is embedded for any $(\t,\a ,\be )$. First note that $M_{\t ,0,0}$ is the toroidal
halfplane layer defined by Karcher in~\cite{ka4}, who proved that $M_{\t ,0,0}$ decomposes in 16 congruent disjoint pieces, each one being the conjugate surface of certain Jenkins-Serrin graph. In particular, $M_{\t ,0,0}$ is embedded
for each $\t \in (0,\frac{\pi }{2})$. Since for $\t $ fixed the heights of the ends of $M_{\t ,\a ,\be }$ depend continuously on $(\a ,\be )$ in the connected set $[0,\frac{\pi }{2}]^2-\{ (0,\t )\} $ (this is clear in the $\xi $-plane model), a standard application of the maximum principle insures that $M_{\t ,\a ,\be }$ is embedded for all values of $\t ,\a ,\be $.

\begin{remark}
Since $P_A$ does not necessarily point to the $x_2$-axis, we must possibly rotate $M_{\t ,\a ,\be }$ by a suitable angle around the $x_3$-axis in order to see it in $\cal S$.
\end{remark}
\par
\noindent
Next we study the conjugate surface $M_{\t ,\a ,\be }^*$ of $M_{\t ,\a ,\be }$. Since
the flux (resp. the period) of the conjugate surface along a given curve in the parameter domain equals the period (resp. the opposite of the flux) of the original surface along the same curve, we deduce from
(\ref{eq:periodendsnosim}), (\ref{eq:fluxg1Mtab}) and (\ref{eq:Perg2Mtab}) that
the period of $M_{\t ,\a ,\be }^*$ at its ends is $\pm \pi \mu (E(\t ,\a ,\be ),0)\in \C \times \R $,
its period along $\g _1$ is $\pi \mu (E(\t ,\a ,\be ),0)$ and the third component of the period of $M_{\t ,\a ,\be }^*$
along $\g _2$ is $-2\pi$. Therefore, $M_{\t ,\a ,\be }^*$ is a complete immersed doubly periodic torus with four horizontal embedded Scherk-type ends, whose period lattice is generated by $\pi \mu (E(\t ,\a ,0),0)$ and $(F(\g _2),2\pi)$, where $F(\g _2)\in \C $ denotes the horizontal part of the flux of $M_{\t ,\a ,\be }$ along $\g _2$. Again any $M_{\t ,\a ,\be }^*$ is embedded by the maximum principle and the embeddedness of $M_{\t ,0,0}^*$ (which decomposes in 16 congruent Jenkins-Serrin graphs). Also note that in order to see $M_{\t ,\a ,\be }^*$ inside the normalized space ${\cal S}$, we must rotate this surface around the $x_3$-axis and rescale it suitably. After this identification, the curve with period zero in the sense of item~{\it 6} of Proposition~\ref{propos1} can be taken as $\g_2^*=\g_1+\g_A\subset M_{\t ,\a ,\be }^*$.

We have defined two families of examples ${\cal A}=\{ M_{\t ,\a ,\be }\} $, ${\cal A}^*=\{ M_{\t ,\a ,\be }^*\}$ inside ${\cal S}$, with $(\t ,\a ,\be )$ varying in ${\cal I}=\left\{(\t,\a,\be)\in (0,\frac{\pi}{2})\times [0,\frac{\pi }{2}]^2\ | \ (\a ,\be )\neq (0,\t )\right\} $.
Clearly this definition can be extended to larger ranges in $(\t,\a,\be)$, but such an extension only produces symmetric images of these surfaces with respect to certain planes orthogonal to the $x_1,x_2$ or $x_3$-axes.
Nevertheless, some of these geometrically equivalent surfaces are considered as distinct points in the space $\widetilde{\cal S}$ defined in Section~\ref{secprelim}.
Another interesting property is that ${\cal A}^*={\cal A}$ (modulo the aforementioned identifications), as the following Lemma shows. In this sense, we can assure that the space of standard examples is {\it self-conjugate.}
\begin{lemma}
\label{lemaautoconj}
Given $(\t ,\a ,\be )\in {\cal I}$, the surface $M_{\t ,\a ,\be }^*$ coincides with $M_{\frac{\pi }{2}-\t ,\a ,\frac{\pi }{2}-\be }$ up to a symmetry in a plane orthogonal to the $x_2$-axis.
\end{lemma}
\begin{proof}
Fix $(\t ,\a ,\be )\in{\cal I}$.
By direct calculation, $\Sigma_{\frac{\pi }{2}-\t }=\{ (\widetilde{z},\widetilde{w})\ | \ \widetilde{w}^2=
(\widetilde{z}^2-1)^2+4\widetilde{z}^2\sec^2\t \} $. Since the Möbius transformation $\varphi (z)=\frac{1-iz}{-i+z}$
applies the set of branch points of the $z$-projection of $\Sigma_\t $ bijectively on the set of branch points
of the $\widetilde{z}$-projection of $\Sigma_{\frac{\pi }{2}-\t }$, it follows that $\Delta (z,w)=(\varphi (z),\widetilde{w}(\varphi (z)))$ is a biholomorphism between $\Sigma_\t $ and $\Sigma_{\frac{\pi }{2}-\t }$.
A simple thought of the role of the angle $\be $ gives that if we extend the definition of the standard examples
to $\be '\in [-\frac{\pi }{2},0)$, then $M_{\t ,\a ,\be '}$ is nothing but the reflected image of
$M_{\t ,\a ,-\be '}$
with respect to the $(x_1,x_3)$-plane (up to a translation). Furthermore, it is straightforward to check that
$g_{\t ,\a ,-\be }=g_{\frac{\pi }{2}-\t ,\a ,\frac{\pi }{2}-\be }\circ \Delta $, where the subindex for $g$ means the parameter angles of the standard example whose Gauss map is $g$. Denoting by $dh_{\t }$ its height differential
(which only depends on $\t $),
a direct computation gives $\Delta^*dh_{\frac{\pi }{2}-\t }=\mu (\frac{\pi }{2}-\t )\frac{\varphi '(z)\, dz}{\widetilde{w}(\varphi (z))}=i\frac{\mu (\frac{\pi}{2}-\t )}{\tan \t }\frac{dz}{w}=i\frac{\mu (\frac{\pi }{2}-\t )}{\mu (\t )\tan \t }\, dh_{\t }$.
Now the Lemma is proved.
\end{proof}

One could ask for which angles $(\t,\a,\be)$ the surface $M_{\t,\a,\be}$ is self-conjugate (i.e. when it is congruent to its conjugate surface).
First note that for any $\t \in (0,\frac{\pi }{2})$, $M_{\t ,\frac{\pi }{2},\be }$ does not depend on $\be \in [0,\frac{\pi }{2}]$ since its spherical configuration is just the rotation image of the one of $M_{\t ,\frac{\pi }{2},0}$ around the $x_3$-axis by angle $\be$.
Having this in mind and using Lemma~\ref{lemaautoconj} we have that the range of angles for which $M_{\t,\a,\be}$ is self-conjugate is $\left(\{\frac{\pi}{4}\}\times (0,\frac{\pi}{2}]\times\{\frac{\pi}{4}\}\right)\cup\left(\{\frac{\pi}{4}\}\times\{\frac{\pi}{2}\}\times[0,\frac{\pi}{2}]\right)$.

As the branch values of the Gauss map of any $M_{\t ,\a ,\be }$ lie on a spherical equator, a result by Montiel and Ros~\cite{mro1} insures that the space of bounded Jacobi functions on $M_{\t ,\a ,\be }$ is $3$-dimensional (they reduce to the linear functions of the Gauss map), a condition usually referred in literature as the {\it nondegeneracy} of $M_{\t ,\a ,\be }$. This nondegeneracy can be interpreted by means of an Implicit Function Theorem argument to obtain that around $M_{\t ,\a ,\be}$, the space ${\cal S}$ is a $3$-dimensional real analytic manifold (Hauswirth and
Traizet~\cite{HausTraizet1}); in particular, the only elements in ${\cal S}$ around a standard example are themselves standard. This local uniqueness result will be extended in the large by Theorem~\ref{thm1} in this paper.

\begin{figure}
\centerline{\includegraphics[width=15.2cm,height=4.1cm]{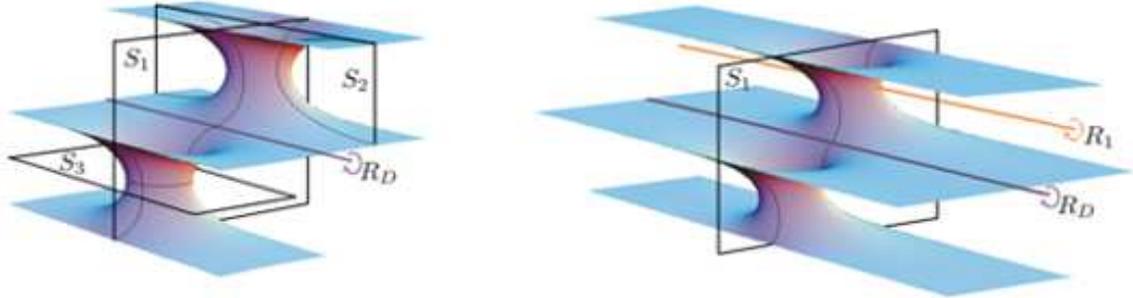}}
\caption{Left: The toroidal halfplane layer $M_{\pi/4,0,0}$. Right: The surface
$M_{\pi/4 ,0,\pi/8}$.}
\label{StandardExamples}
\end{figure}

\begin{figure}
\centerline{\includegraphics[width=13cm,height=4.6cm]{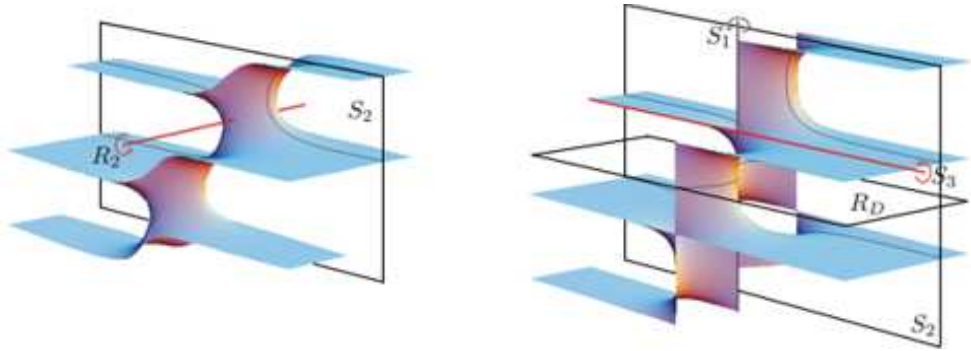}}
\caption{Left: The surface $M_{\pi/4,\pi/4,0}$. Right: $M_{\pi/4 ,\pi /2,0}$, all
whose horizontal level curves are open arcs in $\R^3$.}
\label{StandardExamples2}
\end{figure}

We finish this section by summarizing some additional properties of the standard examples, that can be checked using their Weierstrass representation. For details, see~\cite{PeRo1}.
\begin{proposition}
\label{proposstandardPeRo1}
\begin{description}
\item[]
\item[1.] For any $\t \in (0,\frac{\pi }{2})$, $M_{\t ,0,0}$ admits 3 reflection symmetries
$S_1,S_2,S_3$ in orthogonal planes and contains a straight line parallel to the $x_1$-axis, that induces a $180^{\circ }$-rotation symmetry $R_D$. The isometry group {\rm Iso}$(\t ,0,0)$ is isomorphic to $(\Z /2\Z )^4$, with generators $S_1,S_2,S_3,R_D$ (see Figure~\ref{StandardExamples} left).
\item[2.]  For any $(\t ,\a ) \in (0,\frac{\pi }{2})^2$, $M_{\t ,\a,0}$ is invariant by a
reflection $S_2$ in a plane orthogonal to the $x_2$-axis and by a $180^{\circ }$-rotation $R_2$ around a line parallel to the $x_2$-axis that cuts the surface orthogonally (Figure~\ref{StandardExamples2} left).
{\rm Iso}$(\t ,\a ,0)$ is isomorphic to $(\Z /2\Z )^3$, with
 generators $S_2,R_2, {\cal D}$. Furthermore, {\rm Iso}$(\t ,\frac{\pi }{2},0)={\rm Iso}(\t ,0,0)$ (Figure~\ref{StandardExamples2} right).
\item[3.]  For any $(\t ,\be ) \in (0,\frac{\pi }{2})^2-\{ (\t ,\t )\} $, $M_{\t ,0,\be}$ is invariant by a reflection $S_1$ in a plane orthogonal to the $x_1$-axis, by a $180^{\circ }$-rotation symmetry $R_D$ around a straight line parallel to the $x_1$-axis contained in the surface, and by a $180^{\circ }$-rotation $R_1$ around another line parallel to the $x_1$-axis that cuts the surface orthogonally. {\rm Iso}$(\t ,0,\be )$ is isomorphic to $(\Z /2\Z )^3$, with generators $S_1,R_D,R_1$ (Figure~\ref{StandardExamples} right). Furthermore, {\rm Iso}$(\t ,0,\frac{\pi }{2})=
${\rm Iso}$(\t ,0,0)$.
\item[4.] For any $(\t ,\a ,\be )\in (0,\frac{\pi }{2})^3$, {\rm Iso}$(\t ,\a ,\be )$ is isomorphic to
$(\Z /2\Z )^2$ with generators ${\cal D},{\cal E}$.
\item[5.] For any $(\t ,\a )\in (0,\frac{\pi }{2})^2$, $M_{\t ,\a ,\frac{\pi }{2}}$ is invariant by a reflection $S_3$
in a plane orthogonal to the $x_2$-axis and by a $180^{\circ }$-rotation $R_3$ around a line parallel to the $x_2$-axis that intersects $M_{\t ,\a ,\frac{\pi }{2}}$ orthogonally. {\rm Iso}$(\t ,\a ,\frac{\pi }{2})$ is isomorphic to $(\Z /2\Z )^3$, with generators $S_3,R_3,{\cal D}$.
\item[6.] When $(\t ,\a ,\be )\to (0,0,0)$, $M_{\t ,\a,\be }$ converges smoothly to two vertical catenoids, both with flux $(0,0,2\pi )$.
\item[7.] Let $\t _0\in (0,\frac{\pi }{2})$. When $(\t ,\a ,\be )\to (\t _0,0,\t _0)$, $M_{\t ,\a ,\be }$ converges to a
Riemann minimal example with two horizontal ends, vertical part of its flux $2\pi $ and branch values of its Gauss map
at $0,\infty , i \tan \t _0, -i\cot \t _0\in \overline{\C }$.
\item[8.] When $(\t ,\a ,\be )\to (\frac{\pi }{2},0,\frac{\pi }{2})$, $M_{\t ,\a ,\be }$ converges (after blowing up) to two vertical helicoids.
\item[9.] Let $(\a _0,\be _0)\in [\t ,\frac{\pi }{2}]^2-\{ (0,0)\} $. When $(\t ,\a ,\be )\to (0,\a _0,\be _0)$,
$M_{\t ,\a ,0}$ converges to two singly periodic Scherk minimal surfaces, each one with two horizontal ends and two ends forming angle
$\arccos (\cos \a _0\cos \be  _0)$ with the horizontal.
\item[10.] Let $(\a _0,\be _0)\in [\t ,\frac{\pi }{2}]^2-\{ (0,\frac{\pi}{2})\}$. When $(\t ,\a ,\be )\to (\frac{\pi }{2},\a _0,\be _0)$,
$M_{\t ,\a ,0}$ converges (after blowing up) to two doubly periodic Scherk minimal surfaces, each one with two horizontal ends and two ends forming angle $\arccos (\cos \a _0\sin \be  _0)$ with the horizontal.
\end{description}
\end{proposition}

\section{The moduli space $\cal W$ of Weierstrass representations.}
\label{secW}
Any surface in the hypotheses of Theorem~\ref{thm1} can be seen as an element of a finite dimensional complex manifold by
means of the {\it Hurwitz schemes}, a process which endows families of meromorphic functions of prescribed degree with structures of finite dimensional complex manifolds, see for instance~\cite{Fulton}. This general construction is simpler in our setting, as we now explain.

We will call $\cal W$ to  the space of lists $(\M, g, p_1,\ldots,
p_{2k},q_1,\ldots,q_{2k},[\g])$, where $g:\M\to \overline{\C }$ is a meromorphic degree
$2k$ function defined on a torus $\M$, which is unbranched at its zeros
$\{ p_1,\ldots, p_{2k}\} =g^{-1}(0)$ and poles $\{ q_1,\ldots,q_{2k}\} =g^{-1}(\infty )$, and $[\g]$ is a homology
class in $g^{-1}(\C ^*)$ with $[\g ]\neq 0$ in  $H_1(\M,\Z)$. Note that the same map $g$ can be viewed as an
infinite subset of $\cal W$ by considering different orderings on the points
$p_j,q_j$ and different homology classes $[\g]$. This infinite subset associated to
the same $g$ will be discrete with the topology to be defined below. We will
simply denote by $g$ the elements of $\cal W$, which will be
referred to as {\it marked meromorphic maps}.

Next we endow $\cal W$ with a topology. Given $g\in\cal W$, let $b_1,\dots, b_l\in\C^*$, $l\leq 4k$, be the distinct branch
values of $g$ (by the Riemann-Hurwitz formula, $g$ has  $4k$ branch points counting multiplicity). Take $l$ small
pairwise disjoint disks $D_i$ in $\C^*$ centered at the points~$b_i$ and let
$\Omega=\overline{\C }-\cup_{i=1}^l D_i$. By taking the radii of these
disks small enough, we can assume that both $\Omega$ and $g^{-1}(\Omega)$ are
connected and each component of $g^{-1}(D_i)$ is a disk ($1\leq i\leq l$).
Therefore, $g:g^{-1}(\Omega)\to\Omega$ is an unbranched $2k$-sheeted
covering with connected total space, and each component of $g^{-1}(D_i)$ contains at most one branch point of $g$, possibly of high
multiplicity. With these  data $g,D_1,\dots, D_l$ we define a neighborhood $U(g)$ of
$g$ in $\cal W$ as the set of marked meromorphic maps $f\in\cal W$ such that:
\begin{description}
\item [(1)] {\it $f|_{f^{-1}(\Omega)}$ is unbranched and the restrictions
$f:f^{-1}(\Omega)\to\Omega$, $g:g^{-1}(\Omega)\to\Omega$ are isomorphic as covering
maps.} This allows us to identify conformally $f^{-1}(\Omega)$ and $g^{-1}(\Omega)$,
so that zeros (resp. poles) of $g$ identify with zeros
(resp. poles) of $f$.

\item[(2)] {\it Each component of $f^{-1}(D_i)$ is a disk in the torus $f^{-1}(\overline{\C })$.}
These disks are in bijective correspondence with those of $g$ by means of the
identification in {\bf (1)}, and the total branching order of $g,f$ at corresponding
disks necessarily coincides.

\item[(3)] The identification in {\bf (1)} gives a bijection between the zeros (resp. poles) of
$g$ and those of $f$. {\it We impose that the ordering in the set of zeros and poles
of the marked meromorphic map $f$ coincides
 with the corresponding one for $g$ through this bijection.}

\item[(4)] The homology class $[\g]\in H_1(g^{-1}(\C^*),\Z)$ can be represented
by an embedded cycle $\g\in g^{-1}(\Omega)$, which by {\bf (1)} can be also seen as an embedded cycle
 in $f^{-1}(\Omega)$. Such a cycle induces an homology class in
$f^{-1}(\C^*)$, which we denote $[\g]$ as well. It is easy to show that $[\g]\neq 0$
in $H_1(f^{-1}(\C^*),\Z)$. {\it We impose that the homology class associated to the
marked meromorphic map $f$ coincides with $[\g]$.}
\end{description}
\par
\noindent
Next we endow $\cal W$ with a structure of $4k$-dimensional complex manifold. The
topology defined above implies that if $g\in\cal W$ has $4k$
distinct branch values, then every $f\in\cal W$ close enough to $g$ has also $4k$
distinct branch values. In this case, the map that applies this $f$ to its list of
branch values (choosing an ordering) is a local chart for $\cal W$ around
$g$. Around a generic $g\in\cal W$, a local chart can be obtained by exchanging the
branch values of $g$ by the list $(\sigma_1(g),\dots,\sigma_{4k}(g))$, where
$\sigma_i(g)$ is the value of the symmetric elementary polynomial of degree $i$ on
the unordered list of $4k$ (not necessarily distinct) branch values of $g$,
$1\leq i \leq 4k$. These symmetric
elementary polynomials can be considered as globally defined holomorphic functions
$\sigma_i:\cal W\to\C$, $1\leq i\leq 4k$. Also, the map
$(\sigma_1,\dots,\sigma_{4k}):{\cal W}\to \C^{4k}$ is a local diffeomorphism,
hence $\cal W$ can be seen as an open submanifold of $\C^{4k}$.

The following result deals with compact analytic subvarieties of ${\cal W}$. A
subset $V$ of a complex manifold $N$ is said to be an {\it analytic subvariety} if
for any $p\in N$ there exists a neighborhood $U$ of $p$ in $N$ and a finite number
of holomorphic functions $f_1,\ldots ,f_r$ on $U$ such that $U\cap V=\{q\in U \ | \
f_i(q)=0,\ 1\leq i\leq r\} $. The proof of the following lemma can be found
in~\cite{mpr1}, Lemma 4 (although the definition of ${\cal W}$ in~\cite{mpr1} is
different, the proof still works in our case).

\begin{lemma}
\label{lemacptanalsub}
 The only compact analytic subvarieties of ${\cal W}$ are
finite subsets.
\end{lemma}

\subsection{The height differential associated to a marked meromorphic map.}

Consider a marked meromorphic map $g=(\M, g, p_1,\ldots,p_{2k},q_1,\ldots,q_{2k},[\g])\in {\cal W}$. Since the complex
space of holomorphic differentials on $\M$ is $1$-dimensional, there exists a unique holomorphic $1$-form $\phi=\phi(g)$
on $\M$  such that
\begin{equation}
\label{eq:defphi} \int_\g \phi =2\pi i.
\end{equation}
The pair $(g,\phi)$ must be seen as the Weierstrass data of a potential surface in
the setting  of Theorem~\ref{thm1}, defined on $M=g^{-1}(\C^*)$, with Gauss map $g$ and height
differential $\phi$. Equation~(\ref{eq:defphi}) means that
the period of $(g,\phi )$ along $\g$ is horizontal and its flux along $\g$ has
third coordinate $2\pi $.
We will say that $g\in {\cal W}$ {\it closes periods} when there exists $a\in\R^*$ such that (\ref{eq:gcloseends}) and
the first equation in (\ref{eq:gcloseshomology}) hold with $dh=\phi $ and $\g _2=\g $
(note that the second equation in (\ref{eq:gcloseshomology})
holds by definition of $\phi $).

\begin{lemma}
\label{lemagclosesperiods}
If $g\in {\cal W}$ closes periods, then $(g,\phi )$ is the Weierstrass
pair of a properly immersed minimal surface $M\subset \T\times\R$ for a certain flat torus $\T$, with total curvature $8k\pi$ and $4k$ horizontal Scherk-type ends.
Furthermore, the fluxes at the ends $p_1,q_1,\ldots,p_k,q_k$ are equal to $(\pi a,0,0)$ and opposite to the fluxes at $p_{k+1},q_{k+1},\ldots,p_{2k},q_{2k}$
(here $a\in\R^*$ comes from equation {\rm (\ref{eq:gcloseends})}).
\end{lemma}
\begin{proof}
Since $g$ closes periods, it follows that the period vectors of $(g,\phi)$ at the ends are $\pm(0,\pi a,0)$. Hence it suffices to check that if $\a \subset g^{-1}(\C^*)$ is a closed curve with $\{[\a],[\g]\}$ forming a basis of $H_1(g^{-1}(\overline{\C }),\Z )$, then the period vector of $(g,\phi )$ along $\a$ is linearly independent of $(0,\pi a,0)$. This  holds because $\phi$ is a nonzero holomorphic differential on the torus $g^{-1}(\overline{\C })$ with purely imaginary integral along $\g $, hence $\Re \int_{\a }\phi \neq 0$.
\end{proof}

Equation (\ref{eq:gcloseends}) gives $4k$ residue equations on a marked meromorphic map $g\in {\cal W}$ which are necessary for $g$ to close periods with $dh =\phi $. The fact that the sum of the residues of a meromorphic differential on a compact Riemann surface equals zero implies that it suffices to impose (\ref{eq:gcloseends}) for $1\leq j\leq 2k-1$, so we end up with $4k-2$ equations. Provided that the first equation in (\ref{eq:gcloseshomology}) also holds for $g$ with $\g _2=\g $, the horizontal component of the flux of the corresponding immersed minimal surface $M$ in Lemma~\ref{lemagclosesperiods} is given by
\begin{equation}
\label{eq:horizontalflux}
 F(\gamma )=i\int_\g g\phi\in \C\equiv \R^2.
\end{equation}

\subsection{The ligature map.}
\label{subsecligaturemap}
 We define the {\it ligature map} $L:{\cal W}\to \C^{4k}$
as the map that associates to each marked meromorphic map $g\in\cal W$ the
$4k$-tuple
\[
L(g)=\left(\mbox{Res}_{p_1}\left( g^{-1}\phi\right) ,\dots,\mbox{Res}_{p_{2k-1}}\left( g^{-1}\phi\right) ,
\mbox{Res}_{q_1}(g\phi),\dots,\mbox{Res}_{q_{2k-1}}(g\phi),
\int_\g g^{-1}\phi,\int_\g g\phi\right).
\]
As $\phi$ depends holomorphically on $g$ and all the components of $L$ can be
computed as integrals along curves contained in $g^{-1}(\Omega)$ (see the definition
of the topology of $\cal W$), we conclude that $L$ is holomorphic.
We consider the subset of marked meromorphic maps that close periods,
${\cal M}=\{g\in{\cal W}\ | \ L(g)= L_{(a,b)},\ \mbox{with } a\in\R^*,\ b\in\C\}$, where
\[
L_{(a,b)}=\Big(\underbrace{a,\dots,a}_{1\leq j\leq
k},\underbrace{-a,\dots,-a}_{k+1\leq j\leq 2k-1},\underbrace{-a,\dots,-a}_{2k\leq j\leq 3k-1},\underbrace{a,\dots,a}_{3k\leq j\leq 4k-2},b,\overline b\Big)\in\C^{4k}.
\]
The following map $J$ defines a canonical injection from  the topological space $\widetilde{\cal S}$ of marked surfaces into ${\cal M}$,
\[
(M,p_1,\ldots ,p_{2k},q_1,\ldots ,q_{2k},[\g _2])\mapsto J(M)=(g^{-1}(\overline{\C }),g,p_1,\ldots ,p_{2k}, q_1,\ldots
,q_{2k},[\g_2]),
\]
with $g$ being the Gauss map of $M$.

\begin{lemma}
\label{lemmaJ}
 $J:\widetilde{\cal S}\to {\cal M}$ is an embedding, where ${\cal M}$ has the
restricted topology from the one of ${\cal W}$.
Moreover if we identify $\widetilde{\cal S}$ with $J(\widetilde{\cal S})$, then $\widetilde{\cal S}$ is open and closed in ${\cal M}$.
\end{lemma}
\begin{proof}
First consider a sequence $\{M_n\}_n\subset \widetilde{\cal S}$ converging to a marked surface $M\in\widetilde{\cal S}$.
It is a standard fact that the Gauss map $g_n$ of $M_n$ converges uniformly as $n\to \infty$ to the Gauss map $g$ of $M$ on the compactified torus. From here it is not difficult to check that $J(M_n)$ lies in an arbitrarily small neighborhood of $J(M)$ in the topology of $\cal W$ for $n$ large enough, so $J$ is continuous.
Both the continuity of the inverse map $J^{-1}$ and the openness
and closeness of $\widetilde{\cal S}$ in ${\cal M}$ can be deduced from an argument similar to
the one in the proof of Lemma $6$ of~\cite{mpr1} (in fact, our setting is easier
than~\cite{mpr1} because the Gauss map of any surface in ${\cal S}$ is
unbranched at the ends).
\end{proof}

In the sequel, we will identify $J$ with the inclusion map and see $\widetilde{\cal S}$ as a subset of ${\cal M}$.

\begin{proposition}
\label{proposanalsubv}
Let $H\in \R^3-\{ 0\} $ be a vector parallel to the $x_2$-axis and $F\in \C $. Then,
the set of marked surfaces $M\in \widetilde{\cal S}$ whose periods at the ends are  $\pm  H$ and whose
flux vector along the homology class in the last component of $M$ equals $(F,2\pi )\in \C \times \R \equiv \R^3$ is an analytic subvariety of ${\cal W}$.
\end{proposition}
\begin{proof}
Fix $a\in\R^*$ and $b\in \C $.  It follows from (\ref{eq:periodfluxends}) and (\ref{eq:horizontalflux}) that the set ${\cal M}(a,b)=\{g\in {\cal W}\ | \ L(g)=L_{(a,b)}\}$ coincides with the set of immersed minimal surfaces $M\in {\cal M}$ whose periods and fluxes at the ends are
$P_{p_j}=-P_{p_{k+j}}=-P_{q_j}=P_{q_{k+j}}=-(0,\pi a,0)$,
$F_{p_j}=-F_{p_{k+j}}=F_{q_j}=-F_{q_{k+j}}=(\pi a,0,0)$, for $j=1,\ldots ,k$,
and whose flux along the last component $[\g ]$ of the marked meromorphic map $g$ is
$F_\g=(i\overline{b},2\pi )$.
As $L$ is holomorphic, ${\cal M}(a,b)$ is an analytic subvariety of ${\cal W}$. As a
simultaneously open and closed subset of an analytic subvariety is also an analytic
subvariety, Lemma~\ref{lemmaJ} implies that $\widetilde{\cal S}\cap {\cal M}(a,b)$ is an
analytic subvariety of ${\cal W}$, which proves the Proposition.
\end{proof}

\begin{definition}
{\rm The value of the ligature map $L$
at a marked surface $M\in \widetilde{\cal S}$ is determined by two numbers $a\in \R^*,b\in \C$ so that $\mbox{\rm Res}_{p_1}\left( g^{-1}\, dh\right) =a$ and $F_{\g_2}=(i\overline{b},2\pi )$.
We define the {\it classifying map} $C:\widetilde{\cal S}\to \R^*\times \C$ by $C(M)=(a,b)$.
}
\end{definition}

\begin{remark}
Let $M\in {\cal S}$ be a geometric surface, seen as two marked surfaces $M_1,M_2\in \widetilde{\cal S}$ with
associated homology classes $[\g _2(M_1)],[\g _2(M_2)]\in H_1(M,\Z )$ such that $[\g _2(M_1)]=[\g _2(M_2)]$ in $H_1(\M ,\Z )$ (here $\M$ is the compactification of $M$).
Then, $\g _2(M_1)\cup \g _2(M_2)$ bounds an even number of ends whose periods add up to zero, and the components of $C$ at $M_1,M_2$ satisfy $a(M_1)=\pm a(M_2)\in \R^*$, $b(M_1)=b(M_2)+m\pi a(M_1)$ with $m\in \Z $ even.
\end{remark}

\section{Properness.}
\label{seccurvestim}
In the sequel, we will denote by $K_\Sigma$ the Gaussian
curvature function of any surface $\Sigma$.

\begin{lemma}\label{entornotub}
Let $\cal P$ be a rank $2$ lattice in $\R^3$ and $\Sigma\subset \R^3/\cal P$ a
properly embedded nonflat (orientable) minimal surface with finite topology and $4k$
horizontal Scherk-type ends. Denote by $H, T$ two generators of ${\cal P}$ and assume
that $H$ points to the $x_2$-axis. Suppose also that $|K_\Sigma|\leq c$
where $c>0$. Then,
\begin{description}
\item[{\it (i)}] Both $\| H\|$ and the vertical distance between consecutive left (resp. right) ends of $\Sigma$ are not less than $2/\sqrt{c}$, and the third coordinate of $T$ satisfies
$| T_3|\geq 4k/\sqrt{c}$.
\item[{\it (ii)}] The injectivity radius of $\R^3/{\cal P}$ is bounded by below by $1/\sqrt{c}$ and $\Sigma $
admits a regular open neighborhood of radius $1/\sqrt{c}$ in $\R^3/{\cal P}$.
\end{description}
\end{lemma}
\begin{proof}
Let $\widetilde\Sigma\subset \R^3$ be the connected properly embedded doubly
periodic minimal surface obtained by lifting $\Sigma$.  As $|K_{\widetilde{\Sigma}}|\leq c$
and the ends of $\widetilde{\Sigma }$ are asymptotic to horizontal halfplanes, a standard application
of the maximum principle shows that
$\widetilde \Sigma$ has a regular neighborhood $\widetilde \Sigma(\ve )$ of positive
radius $\ve =1/\sqrt c$ (see~\cite{mr1,ros5} for similar arguments). In particular,
the vertical separation between two consecutive left ends of $\widetilde{\Sigma}$ in
the same fundamental domain is greater than or equal to $2\ve $. As $\Sigma $
has $2k$ left ends, it follows that $| T_3|\geq 4k\ve $. Now take a point
$p\in \widetilde{\Sigma}$ where the normal vector to $\widetilde{\Sigma}$ points
to $(0,1,0)$. $\widetilde{\Sigma}$ can be expressed locally around $p$ as a graph
${\cal G}$ over a disk in the tangent space at $p$ and the same holds around the homologous
point $p+ H$. Reasoning as above with ${\cal G}$, ${\cal G}+H$ instead of consecutive left ends of
$\widetilde{\Sigma }$, we have that $\| H\|\geq 2\ve $. This proves {\it (i)}.

Let $r$ be the injectivity radius of $\R^3/{\cal P}$. Since $\R^3/{\cal P}$ is flat, there are no conjugate
points on geodesics in $\R^3/{\cal P}$ and thus $r$ is half of the minimum length of a closed geodesic in
$\R^3/{\cal P}$. Therefore,
\[
2r=\min_{v\in {\cal P}-\{0\}}\| v\| \geq \min (\|  H\| ,| T_3|)\geq \min \left(
\frac{2}{\sqrt{c}},\frac{4k}{\sqrt{c}}\right) =\frac{2}{\sqrt{c}}.
\]
Finally, any two distinct points $p,q\in\R^3$ homologous by $\cal P$ are separated by a
distance greater than or equal to $2 r$. It follows that $\widetilde\Sigma(\ve)/\cal P$
is a regular neighborhood of $\Sigma$ of radius $1/\sqrt{c}$.
\end{proof}

Next we describe all possible limits of sequences $\{M_n\}_n$ of surfaces
under the hypotheses of Theorem~\ref{thm1},
with the additional assumption of having
uniform curvature bounds. For such a surface $M_n$, Proposition~\ref{propos1}
insures that its period lattice ${\cal P}_n$ is generated by the horizontal period vector $ H_n=(0,\pi a_n,0)$ at the ends of $M_n$ and by a nonhorizontal vector
 $T_n=P_{\g_1(n)}\neq 0$, where $a_n>0$, $\g _1(n)\subset M_n$ is a closed curve so that
$[\g _1(n)]\neq 0$ in $H_1(\M_n,\Z )$ and $\M _n$ is the compactification of $M_n$.
We will use this notation throughout the proof of the following Proposition.

\begin{proposition}
\label{proposdifferentlimits}
 Let $\{\widetilde M_n\}_n$ be a sequence of properly
embedded doubly periodic minimal surfaces, each $\widetilde M_n$ invariant by a rank
$2$ lattice ${\cal P}_n$ as above. Suppose that for all~$n$, $M_n=\widetilde M_n/{\cal P}_n$
has genus $1$ and $4k$ horizontal Scherk-type ends, $\widetilde M_n$ passes
through the origin of $\R^3$ and $|K_{\widetilde M_n}(0)|=1$ is a maximum value of
$|K_{\widetilde M_n}|$. Then (after passing to a subsequence), $\widetilde M_n$
converges uniformly on compact subsets of $\R^3$ with multiplicity $1$ to a properly
embedded minimal surface $\widetilde M_\infty\subset\R^3$ which lies in one of the
following cases:
\begin{description}
\item [{\it (i)}] $\widetilde M_\infty$ is a vertical catenoid with flux $(0,0,2\pi)$. In this case,
 both $\{H_n\}_n,\{T_n\}_n$ are unbounded for any choice of $T_n$ as above.
\item [{\it (ii)}] $\widetilde M_\infty$ is a vertical helicoid with period vector $(0,0,2\pi m)$ for some $m\in\N$.
Now $\{H_n\}_n$ is unbounded and there exists a choice of $T_n$ for which $\{T_n\}_n\to (0,0,2\pi m)$ as $n\to \infty$.
\item [{\it (iii)}] $\widetilde M_\infty$ is a Riemann minimal example with horizontal
ends. Moreover, $\{H_n\}_n$ is unbounded and certain choice of $\{T_n\}_n$ converges to the period vector of $\widetilde M_\infty$.
\item [{\it (iv)}] $\widetilde M_\infty$ is a singly periodic Scherk minimal surface, two of whose ends
are horizontal. Furthermore, any choice of $\{T_n\}_n$ is unbounded, $\{H_n\}_n$  converges to the period
vector $H_\infty=(0,a,0)$ of $\widetilde M_\infty$ (with $a>0$), and
$\widetilde{M}_\infty/ H_\infty$ has genus zero.
\item [{\it (v)}] $\widetilde M_\infty$ is a doubly periodic Scherk minimal surface. In this case, $\{H_n\}_n$, $\{T_n\}_n$  converge respectively to period vectors $H_{\infty },T_{\infty }$ of $\widetilde M_\infty$, and
$\widetilde{M_{\infty }}/{\{ H_{\infty },T_{\infty }\} }$ has genus zero with at least two
horizontal ends and exactly two nonhorizontal ends.
\item [{\it (vi)}] $\widetilde M_\infty$ is a doubly periodic minimal surface invariant by a
rank~$2$ lattice ${\cal P}_\infty$, $M_\infty=\widetilde M_\infty/{\cal
P}_\infty$ has genus one and $4k$ horizontal Scherk-type ends, and $\{H_n\}_n\to
 H_\infty$, $\{T_n\}_n\to  T_\infty$, where $H_\infty, T_\infty$ are defined by Proposition~\ref{propos1} applied to
$M_\infty$.
\end{description}
\end{proposition}
\begin{proof}
By Lemma~\ref{entornotub}, $\{\widetilde M_n\}_n$ has local uniform area bounds. As
$\{|K_{\widetilde M_n}|\}_n$ is uniformly bounded and the origin is an accumulation
point of $\{\widetilde M_n\}_n$, after choosing a subsequence
$\{\widetilde M_n\}_n$ converges uniformly on compact subsets of $\R^3$ to a
properly embedded minimal surface $\widetilde M_\infty$ with $0\in\widetilde
M_\infty$ and $|K_{\widetilde M_\infty}(0)|=1$. Since $\widetilde M_\infty$ is complete,
orientable and not flat, it is not stable (see e.g. do Carmo and Peng \cite{cp1} or Fischer-Colbrie and Schoen~\cite{fs1}). This
implies that the multiplicity of the limit $\{\widetilde M_n\}_n\to\widetilde
M_\infty$ is one, see for instance \cite{pro2}. We now discuss on the
boundedness of the sequences $\{H_n\}_n,\{T_n\}_n$.

{\bf Suppose that both $\{ H_n\}_n$ and any choice of $\{T_n\}_n$ are unbounded}.
After passing to a subsequence, we can assume that the minimum length of a nonzero vector in ${\cal P}_n$ diverges to $\infty$ as $n\to \infty $.
Thus $\widetilde M_\infty$ can be viewed as limit of surfaces with boundary $\widetilde
\Sigma_n$ contained in fundamental domains of $\widetilde M_n$, so $\widetilde
M_\infty$ has finite total curvature. As each $\widetilde \Sigma_n$ has genus zero,
the same holds for $\widetilde M_\infty$. Therefore $\widetilde M_\infty$ is a
catenoid (López and Ros \cite{lor1}). As the Gauss map of $\widetilde M_n$ omits the
vertical directions for all $n$, the same holds for $\widetilde{M}_{\infty }$ by the Open Mapping Theorem. Hence
$\widetilde M_\infty$ is a vertical catenoid.
Since $|K_{\widetilde M_\infty}|$ reaches a maximum at the origin with value $1$, the flux of $\widetilde M_\infty$
must be $(0,0,2\pi)$.

{\bf Suppose that $\{H_n\}_n$ is unbounded and there is a bounded choice of $\{T_n\}_n$}.
We do not
loss generality assuming that $\{H_n\}_n\to \infty$ and $\{T_n\}_n\to
 T_\infty\in \R^3$. By Lemma~\ref{entornotub}, $|( T_n)_3|\geq 4k$ for all
$n$, thus $|( T_\infty)_3|\geq 4k$. In particular $T_\infty\neq 0$ and
$\widetilde M_\infty$ is singly periodic, invariant by the translation by $T_\infty$. As the Gauss map of $\widetilde
M_\infty/ T_\infty$ has degree at most $2k$ and this surface is properly
embedded in $\R^3/ T_\infty$, it must have finite topology, and all its ends
are simultaneously asymptotic to ends of planes, helicoids or Scherk minimal surfaces.
Note that $\widetilde M_\infty/ T_\infty$ cannot have Scherk-type ends, because otherwise
$\widetilde M_\infty/ T_\infty$ would have at least two horizontal Scherk-type
ends (the Gauss map of $\widetilde M_\infty/ T_\infty$ omits the
vertical directions), so $T_\infty$ would be horizontal, a
contradiction.

Next we prove that if the ends of $\widetilde M_\infty/ T_\infty$ are helicoidal, then
$\widetilde M_\infty$ is a vertical helicoid. As the Gauss map $g_{\infty }$
of $\widetilde M_\infty/ T_\infty$ misses $0,\infty $, we deduce that
the height differential $dh_{\infty }$ of $\widetilde M_\infty/ T_\infty$ has no zeros
and all the ends of $\widetilde M_\infty/ T_\infty$ are asymptotic to vertical helicoids, where $dh_{\infty }$ has simple poles. Calling $P$ to the number of ends of $\widetilde M_\infty/ T_\infty$ and $G$ to the genus of its compactification, it follows that $P+2(G-1)=0$ from where $P=2$ and $G=0$. In this situation, it is well known that
$\widetilde M_\infty$ is a vertical helicoid (Toubiana~\cite{to1} or Pérez and Ros~\cite{pro1}). Furthermore,
its period vector is of the form $(0,0,2\pi m)$ with $m\in \N$ because the
maximum absolute curvature of $\widetilde{M}_\infty/ T_\infty$ is one. Thus we have proved
{\it (ii)}.

Now assume that the ends of $\widetilde{M}_\infty/ T_\infty$ are planar. A direct lifting
argument and the maximum principle insure that $\widetilde{M}_\infty/ T_\infty$ cannot have
genus zero, thus it is a properly embedded minimal torus in $\R^3/ T_\infty$ with
a finite number of planar ends. By a Theorem of Meeks, Pérez and Ros~\cite{mpr1},
$\widetilde{M}_\infty/ T_\infty$ is a Riemann minimal example. As its Gauss map misses
 the vertical directions, the ends of $\widetilde{M}_\infty/ T_\infty$ are
horizontal and we have {\it (iii)}.

{\bf Suppose that $\{H_n\}_n$ is bounded and any choice of $\{T_n\}_n$ is unbounded.}
After extracting a subsequence, we can
assume $\{H_n\}_n\to  H_\infty=(0,a,0)$ with $a\in [0,\infty )$ and $\{
 T_n\}_n\to \infty $. By Lemma~\ref{entornotub}, $a\geq 2$ and
$\widetilde{M}_\infty$ is singly periodic, invariant by translation by $H_\infty$. The
well defined Gauss map $g_{\infty }$ of $\widetilde{M}_\infty / H_\infty$
has degree less than or equal to $2k$, hence $\widetilde{M}_\infty/ H_\infty$
has finite total curvature. As it is embedded and proper, all ends of $\widetilde{M}_\infty / H_\infty$
are again simultaneously planar, helicoidal or of Scherk-type. Moreover,
$\widetilde{M}_\infty/ H_\infty$ has at least two ends with vertical limit
normal vector because $g_{\infty }$ omits $0,\infty $. As $ H_\infty$ is horizontal,
properness of $\widetilde{M}_\infty/ H_\infty$
prevents its ends to be planar. Also, its ends are not helicoidal because in such
case the period vectors at the ends would be vertical. Hence
$\widetilde{M}_\infty/ H_\infty$ has Scherk-type ends
(at least two of them horizontal) and genus zero or one.

\begin{assertion}
\label{asserunbranchedends}
In the above situation, the well defined height differential $dh_{\infty }$ of $\widetilde{M}_\infty/ H_\infty$ cannot have a zero at any horizontal end (equivalently, $g_{\infty }$ is unbranched at the horizontal ends of
$\widetilde{M}_\infty/ H_\infty$).
\end{assertion}

To see the Assertion, suppose on the contrary that $g_{\infty }$ is branched at a puncture $p_{\infty }$ corresponding
to a horizontal Scherk-type end of $\widetilde{M}_\infty/ H_\infty$, say with value zero.
Since the Gauss map $g_n$ of $\widetilde{M}_n/H_n$ is unbranched at its zeros, there exists an integer $l\geq 2$
such that $l$ distinct zeros $p_1(n),\ldots ,p_l(n)$ of $g_n$ converge to $p_{\infty }$ (we can think of pieces of the compactifications of $\widetilde{M}_{\infty }/T_{\infty }, \widetilde{M}_n/T_n$ defined on a common disk $D$ centered at $p_{\infty }$ so that $p_1(n),\ldots ,p_l(n)\in D $ for $n$ large). Moreover, the heights of the ends of $\widetilde{M}_n/H_n$ corresponding to $p_1(n),\ldots ,p_l(n)$ converge to the height of the end $p_{\infty }$ of
$\widetilde{M}_{\infty }/H_{\infty }$. By Lemma~\ref{entornotub}, we deduce that
$l=2$ and $p_1(n),p_2(n)$ are a left end and a right end, and so their fluxes are opposite and the same holds for their periods. Therefore both the period and flux of $\widetilde{M}_{\infty }/H_{\infty }$ at $p_{\infty }$ vanish, a contradiction. This proves Assertion~\ref{asserunbranchedends}.

Note that $dh_\infty$ has a pole at each nonhorizontal Scherk-type end of $\widetilde{M}_\infty/ H_\infty$
(and there are at least two of these
ends by the maximum principle). Using Assertion~\ref{asserunbranchedends} and the fact that $g_{\infty }$ misses
$0,\infty$ on $\widetilde{M}_\infty/ H_\infty$, we conclude that $dh_{\infty }$  has no zeros at the compactification
of $\widetilde{M}_\infty/ H_\infty$. Thus $\widetilde{M}_\infty/ H_\infty$ has genus zero and $dh_{\infty }$ has
exactly two poles or equivalently, $\widetilde{M}_\infty/ H_\infty$ has exactly two nonhorizontal Scherk-type ends.
In this situation, Pérez and Traizet~\cite{PeTra1} have proved that
$\widetilde{M}_\infty/ H_\infty$ is a singly periodic Scherk minimal surface of genus zero. This shows {\it (iv)}.

{\bf Suppose that $\{H_n\}_n$ is bounded and there exists a bounded choice of $\{T_n\}_n$}. As before,
it can be supposed that $\{H_n\}_n\to  H_\infty=(0,a,0)$ with $a\geq 2$ and $\{T_n\} _n\to  T_\infty\in \R^3$
with $|( T_\infty)_3|\geq 4k$. Thus, $\widetilde{M}_\infty$ is doubly periodic, invariant by the rank 2 lattice
${\cal P}_\infty$ generated by $H_\infty, T_\infty$ and $\widetilde{M}_\infty/{\cal P}_\infty$ has genus zero
or one with a finite number of Scherk-type ends.

Since the Gauss map $g_\infty$ of $\widetilde{M}_\infty/{\cal P}_\infty$ misses $0,\infty$, this surface has at least two horizontal ends. A suitable modification of Assertion~\ref{asserunbranchedends} gives that the height differential $dh_\infty$ of $\widetilde{M}_\infty/{\cal P}_\infty$ (which is also well defined) does not vanish at any of the horizontal ends of this surface. Therefore either $\widetilde{M}_\infty/{\cal P}_\infty$ has genus zero with exactly two nonhorizontal Scherk-type ends, or it has genus one and all its ends are horizontal.
In the first case, a theorem by Wei~\cite{wei2} (see also Lazard-Holly and Meeks~\cite{lhm}) insures that $\widetilde{M}_\infty$ is a doubly periodic Scherk minimal surface, so {\it (v)} is proved.
Finally, assume that $\widetilde{M}_\infty/{\cal P}_\infty$ has genus one. Since the
total curvature of $M_n$ is $8k\pi $,
$\widetilde{M}_\infty/{\cal P}_\infty$ must have total curvature at most $8k\pi$.
Using the Meeks-Rosenberg formula~\cite{mr3} and the fact that $\widetilde{M}_\infty/{\cal P}_\infty$ is
properly embedded with genus one and Scherk-type ends, we deduce that
$\widetilde{M}_\infty/{\cal P}_{\infty }$ has at most $4k$ ends.
On the other hand, the spacing between left (resp. right) ends of $M_n$ is bounded
away from zero by Lemma~\ref{entornotub}, which implies that the $2k$ heights
corresponding to the left (resp. right) ends of $M_n$ converge to $2k$ distinct heights in
$\R^3/ H_\infty$.
Since noncompact horizontal level sets of the $\widetilde{M}_n/ H_n$ tend to noncompact level sets of $\widetilde{M}_\infty/ H_\infty$, we conclude that $\widetilde{M}_{\infty}/{\cal P}_\infty$ has exactly $4k$ horizontal ends. Now the proof is complete.
\end{proof}

The following result is a crucial curvature estimate in terms of the classifying map $C$.
\begin{proposition}
\label{proposcurvestim}
 Let $\{M_n\}_n\subset \widetilde{\cal S}$ be a sequence of marked surfaces. Suppose
that  $C(M_n)=(a_n,b_n)\in\R^*\times\C$  satisfies
\begin{description}
    \item[{\it (i)}] $\{ a_n\}_n$ is bounded away from zero.
    \item[{\it (ii)}] $\{ |b_n|\}_n$ is bounded by above.
\end{description}
Then, the sequence of Gaussian curvatures $\{ K_{M_n}\} _n$ is uniformly bounded.
\end{proposition}
\begin{proof}
The proof is based on the one of Theorem 4 in~\cite{mpr1}, so we will only go into the details of what is new in this setting.
By contradiction, assume that $\l
_n:=\max_{M_n}\sqrt{|K_{M_n}|}\to \infty$ as $n\to \infty $. Let $\Sigma_n= \l
_nM_n\subset \R^3/{\l_n\cal P}_n$, where ${\cal P}_n=\mbox{Span}\{H_n, T_n\}$ is
the rank 2 lattice associated to $M_n$ and $H_n$ is the period vector at the ends of $M_n$ (up to sign).
Let us also call $\widehat{\Sigma
}_n,\widetilde{\Sigma}_n$ to the respective liftings of $\Sigma_n$ to $\R^3/\l
_n T_n$ and to $\R^3$.

After translation of $\widetilde{\Sigma}_n$ to have maximum
absolute Gauss curvature one at the origin, Proposition~\ref{proposdifferentlimits} implies that (after passing to a subsequence) $\{\widetilde{\Sigma}_n\}_n$ converges smoothly to a properly embedded minimal surface ${\cal H}_1\subset \R^3$, which must lie in one of the six possibilities in Proposition~\ref{proposdifferentlimits}.
Since $a_n$ is bounded away from zero and $\l_n\to \infty $, the period vectors $\l_n
H_n$ at the ends of $\Sigma_n$ diverge to $\infty $. Therefore ${\cal H}_1$
must be a vertical catenoid, a vertical helicoid or a Riemann minimal example with horizontal ends. If ${\cal H}_1$ were a catenoid or a Riemann minimal example, then the vertical part of the flux of $\widetilde{\Sigma }_n$ along a compact horizontal section, which is $2\pi \l_n$ by item {\it 7} of Proposition~\ref{propos1}, would converge to the vertical part of the flux of ${\cal H}_1$, which is finite.
This contradicts that $\l_n\to \infty$ and so, ${\cal H}_1$ is a
vertical helicoid with period vector $T=(0,0,2\pi m)$ for certain $m\in \N$ (furthermore, we can choose the
period vector $T_n$ of $M_n$ with $\{ \l_nT_n\}_n\to T$ as $n\to \infty $, see Proposition~\ref{proposdifferentlimits}).

Let $\pi_n:\R^3/\l_n {\cal P}_n\rightarrow \{x_3=0\}/\l_n {\cal P}_n$ be the linear
projection  in the direction of $T_n$, $\Pi_n^H:\R^3/\l_nT_n\rightarrow
\R^3/\l_n{\cal P}_n$ the quotient projection modulo $\l_n H_n$ and $\D$ the unit
disk centered at the origin in $\{x_3=0\}$. From now on we will only consider $n$
large so that ${\cal H}_1(n)={\Sigma}_n\cap\pi_n^{-1}(\D/\l_n {\cal P}_n)$ is
connected and extremely close to a piece of ${\cal H}_1$ containing its axis, $\l_n
T_n$ being close to $T$. Let $\overline{\Sigma}_n$ be the torus obtained by
attaching to $\Sigma_n$ its $4k$ ends. Note that ${\cal H}_1(n)$ does not separate
$\overline{\Sigma}_n$, hence ${\cal F}_n=\overline{\Sigma}_n-{\cal H}_1(n)$ is a
compact annulus that contains the $4k$ ends of $\Sigma_n$. Let $N_n$ be the Gauss
map of $\Sigma_n$.

We claim that for all $n$ large, $N_n$ takes horizontal values on ${\cal F}_n$.
To see this, consider the compact intersection $\G_n$ of $\Sigma_n$ with a totally geodesic horizontal cylinder in $\R^3/\l_n {\cal P}_n$ not asymptotic to the ends of $\Sigma_n$.
Viewed in $\overline{\Sigma}_n$, $\G_n$ is a non nulhomotopic embedded closed curve whose intersections with ${\cal H}_1(n),{\cal F}_n$ are two open arcs with common end points $A_n,B_n$. As ${\cal H}_1(n)$ is very close to a piece of a vertical helicoid containing its axis, we deduce that $N_n(A_n),N_n(B_n)$ are in
different hemispheres of $\esf^2$ with respect to the vertical direction. Now our
claim follows by continuity.

Let $\esf^1\subset \esf^2$ be the horizontal equator in the sphere. We next show that for all $\t\in\esf^1$ and $n\in \N$ large enough, ${\cal F}_n\cap N_n^{-1}(\t)$ is not empty.
Since ${\cal H}_1(n)$ is almost a piece of a vertical helicoid containing its axis, $(N_n|_{{\cal H}_1(n)})^{-1}(\esf^1)$ must be a simple closed curve in $\Sigma_n$ (viewed in $\R^3$, $(N_n|_{{\cal H}_1(n)})^{-1}(\esf^1)$ is an open embedded arc whose end points differ in $\l_n T_n$) that covers $\esf^1$ with finite multiplicity through $N_n$.
As $N_n$ is horizontal somewhere in ${\cal F}_n$, we deduce that the number of sheets of this last covering is less than $2k$,
which implies that for all $\t \in \esf^1$, ${\cal F}_n\cap N_n^{-1}(\t )$ is not empty as
desired.

Next we produce a second helicoid as a limit of different translations of the
$\widetilde{\Sigma}_n$.
Since the total branching number of $N_n$ is $4k$, there exists a spherical disk $D_{\esf^2}(\t_n,\ve )\subset \esf^2$ centered at some $\t_n\in\esf^1$ with uniform radius $\ve >0$ which is free of branch values of $N_n$.
By the last paragraph, for $n$ large we can find points $p_n\in \widetilde{\Sigma}_n$ with $N_n(p_n)=\t_n$ and whose projections into $\R^3/\l_n {\cal P}_n$ lie inside ${\cal F}_n$.
The sequence of translated surfaces $\{\widetilde{\Sigma}_n-p_n\}_n$ has uniform
curvature and area bounds and all the $\widetilde{\Sigma}_n-p_n$ have normal vector
$\t_n$ at the origin.
After passing to a subsequence, $\widetilde{\Sigma}_n-p_n$ converges uniformly on compact subsets of $\R^3$ to a properly embedded minimal surface ${\cal H}_2$ having vertical tangent plane at the origin.
${\cal H}_2$ is not flat, because otherwise one could construct an arc $\widetilde\a_n$ contained in $\{x\in \widetilde{\Sigma}_n\ | \ d_{\widetilde{\Sigma}_n}(x,p_n)<\frac{3}{2}\|
T\|\}$, $\widetilde\a_n$ giving rise to the period vector $\l_n T_n$.
By construction, $N_n(\widetilde\a_n)$ would lie in $D_{\esf^2}(\t_n,\ve )$ for $n$
large enough, but this last disk is free of branch values of $N_n$ hence
$N_n^{-1}(D_{\esf^2}(\t_n,\ve ))$ consists of disjoint disks inside $\Sigma_n$,
one of which would contain $\a_n=\widetilde\a_n/\l_n{\cal P}_n$ in contradiction with the fact that $\a_n$ is homotopically nontrivial on $\Sigma_n$. Hence ${\cal H}_2$ is
not flat.
Applying Proposition~\ref{proposdifferentlimits} (suitably modified so
that the absolute Gaussian curvature of the surfaces in the sequence is not 1 at the
origin, but it is uniformly bounded by above and that the limit surface is not a
plane) and our previous arguments to eliminate all limits other than a helicoid, we
conclude that ${\cal H}_2$ is another vertical helicoid. Since both ${\cal
H}_1,{\cal H}_2$ are limits of translations of the $\widetilde{\Sigma}_n$, the
period vector of ${\cal H}_2$ is again $T=\lim_n\l_n T_n$.

Consider two disjoint round disks ${\cal D}_1(n),{\cal D}_2(n)$ in the totally geodesic cylinder $\{x_3=0\}/\l_n{\cal P}_n\subset\R^3/\l_n{\cal P}_n$, with common radius $r_n$
such that the annular component ${\cal H}_i(n)=\Sigma_n\cap\pi_n^{-1}({\cal D}_i(n))$ is arbitrarily close to a translated copy of the forming helicoid ${\cal H}_i/ T$ minus neighborhoods of its ends, $i=1,2$.
After passing to a subsequence, we can also choose $r_n$ so that
\begin{enumerate}
    \item $r_n\to \infty$ as $n\to \infty $.
    \item $\frac{r_n}{\l_n}\to 0$ as $n\to \infty $.
    \item \label{terceracondicion} The normal direction to $\Sigma_n$ along the helix-type curves in the
    boundary of ${\cal H}_i(n)$ makes an angle less than $1/n$ with the vertical, $i=1,2$.
\end{enumerate}
(Note that we have exchanged the former ${\cal H}_1(n)=\Sigma_n\cap\pi_n^{-1}(\D)$ by a bigger one ${\cal H}_1(n)=\Sigma_n\cap \pi_n^{-1}(D_1(n))$ but all the preceding arguments remain valid now).
We claim that the extended Gauss map $N_n$ applies $\overline{\Sigma}_n-({\cal H}_1(n)\cup
{\cal H}_2(n))$ in the spherical disks centered at the North and South Poles of $\esf^2$ with radius $1/n$.
Since $N_n$ is an open map, the condition~\ref{terceracondicion} above shows that it suffices to check that $\Sigma_n-({\cal H}_1(n)\cup {\cal H}_2(n))$ has no points with horizontal normal vector. If horizontal
normal vectors occurred in $\Sigma_n-({\cal H}_1(n)\cup {\cal H}_2(n))$, then our former
arguments would give further vertical helicoids ${\cal H}_3,\ldots ,{\cal H}_s$ as limits of
translations of subsequences of the $\widetilde \Sigma_n$ (finitely many limits
because each one consumes at least $-4\pi $ of total curvature, and the total curvature of the $\Sigma_n$ is fixed
$-8k\pi$).
Therefore we have $s$ pairwise disjoint annuli ${\cal H}_j(n)=\Sigma_n\cap\pi_n^{-1}({\cal D}_j(n))$, each one arbitrarily close to the corresponding vertical helicoid ${\cal H}_j$ minus a neighborhood of its two ends, where each ${\cal D}_j(n)\subset\{x_3=0\}/\l_n {\cal P}_n$ is a round disk as before, $j=1,\ldots,s$.
In this setting, the complement of ${\cal H}_1(n)\cup\ldots \cup {\cal H}_s(n)$ in $\overline{\Sigma}_n$ would consist of $s$ closed annuli ${\cal F}_1(n),\ldots ,{\cal F}_s(n)$ labelled so that
${\cal H}_i(n)$ is consecutive to ${\cal F}_i(n)$ and all ${\cal H}_i(n),{\cal F}_i(n)$ are disposed cyclically
for  $i=1,\ldots ,s$. Moreover, $N_n({\cal F}_1(n)\cup \ldots \cup {\cal F}_s(n))$ omits $\esf^1$
hence $N_n({\cal F}_1(n)\cup \ldots \cup {\cal F}_s(n))$ is contained in a small neighborhood of
the vertical directions in $\esf^2$. When restricted to the surface ${\cal F}_i(n)$, the
projection $\pi_n$ extends smoothly through the ends producing a map
\begin{equation}
\label{eq:coveringnoncompact}
f_n:=\pi_n|_{{\cal F}_i(n)}:{\cal F}_i(n)\rightarrow \frac{\{x_3=0\}\cup \{\infty \}}{\l_n {\cal P}_n}- \cup_{j=1}^s{\cal D}_j(n).
\end{equation}
$f_n$ is a proper local diffeomorphism, hence a finite sheeted covering map.
Since $\partial {\cal F}_i(n)$ has two components it follows that $s=2$, a contradiction.
Thus $N_n[\overline{\Sigma}_n-({\cal H}_1(n)\cup {\cal H}_2(n))]\subset D_{\esf^2}(*,1/n)$ where $*=0,\infty $.

So far, we have proved that the two components ${\cal F}_1(n),{\cal F}_2(n)$ of $\overline{\Sigma}_n- ({\cal H}_1(n)\cup {\cal H}_2(n))$ are closed annuli, each one is noncompact when viewed in $\Sigma_n$.
The number of ends of $\Sigma_n$ in ${\cal F}_i(n)$ is twice the number of sheets $\#_i$ of the covering map $f_n$ in
(\ref{eq:coveringnoncompact}).
Moreover, the boundary components $\a_{i,1}(n),\a_{i,2}(n)$ of ${\cal F}_i(n)$ apply by $f_n$ respectively on the circumferences $\partial {\cal D}_1(n),\partial {\cal D}_2(n)$, both with multiplicity $\#_i$. Since $\a_{1,j}(n),\a_{2,j}(n)$ are the helix-type boundary curves of the forming helicoid ${\cal H}_j(n)$, we deduce that $\#_1=\#_2$.
Since there are the same number of left and right ends of $\Sigma_n$ in ${\cal F}_i(n)$, it follows that $\#_1=k$, see Figure~\ref{magda1}.
\begin{figure}
\centerline{\includegraphics[width=12.7cm,height=5cm]{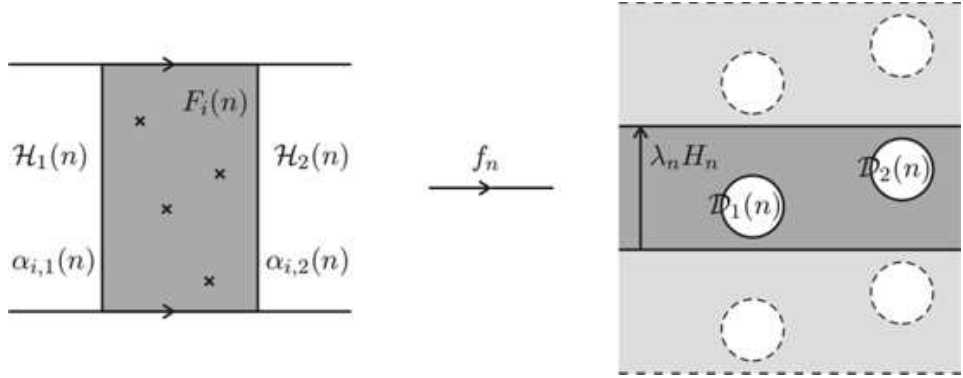}}
\caption{In this
case, the covering map $f_n$ has two sheets.}
\label{magda1}
\end{figure}

In this setting and coming back to the original scale, one can follow the argument in pages 117-118 of~\cite{mpr1} to
construct a closed embedded curve $\G_n\subset M_n$ formed by four consecutive arcs
$L_1(n)^{-1}*\be_1(n)*L_2(n)*\be_2(n)$, where $L_1(n),L_2(n)$ are liftings of the distance minimizing horizontal segment $L(n)$ from $\frac{1}{\l_n}\partial {\cal D}_1(n)$ to
$\frac{1}{\l_n}\partial {\cal D}_2(n)$, lying in consecutive sheets by the covering
\[
f_n:{\textstyle \frac{1}{\l_n}[{\cal F}_1(n)\cup {\cal F}_2(n)]\to \frac{\{x_3=0\}\cup \{\infty \}}{{\cal P}_n}-\frac{1}{\l_n}[{\cal D}_1(n)\cup {\cal D}_2(n)]}
\]
and $\be_1(n),\be_2(n)$ are small arcs contained
in $\frac{1}{\l_n}{\cal H}_1(n),\frac{1}{\l_n}{\cal H}_2(n)$ respectively
(we abuse of notation keeping the label $f_n$ in the original scale
of $M_n$, although it has been defined in the scale of $\Sigma_n=\l_nM_n$).
For later uses, we will describe $\be_1(n),\be_2(n)$ more precisely. $\be_1(n)$ consists of the union of three consecutive arcs $l_1,l_2,l_3$, where $l_1, l_3$ are at almost constant height and $l_2\subset N_n^{-1}(\esf^1)$, and $\be_2(n)$ is similarly defined, see Figure~\ref{helicoids}.

\begin{figure}
\centerline{\includegraphics[width=14cm,height=6.67cm]{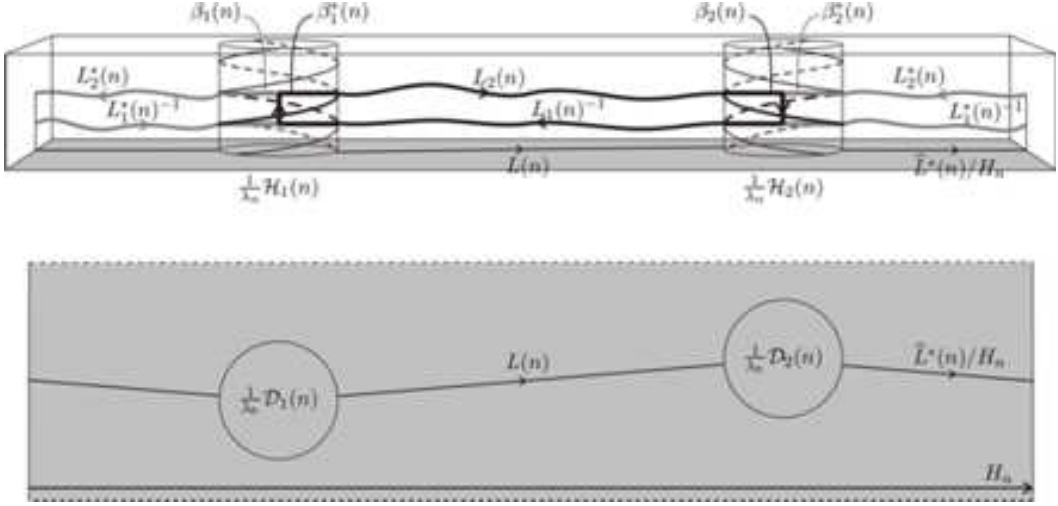}}
\caption{Front and top views of two helicoids forming with the direction of $L(n)$
tending to the direction of the $x_2$-axis. In the front view, the sided ends are
placed ahead and behind the box that contains both helicoids.}
\label{helicoids}
\end{figure}

Next we study the relation between the homology classes $[\G _n],[\g _2(n)]\in H_1(M_n,\Z )$, where
$[\g _2(n)]$ is the last component of the marked surface $M_n\in \widetilde{\cal S}$ (recall that $C(M_n)=(a_n,b_n)$ where $F_{\g_2(n)}=(F(\g _2(n)),2\pi )=(i\overline{b_n},2\pi )$ is the flux vector of $M_n$ along $\g_2(n)$).
Let $g_n$ be the complex Gauss map of $M_n$.
Since $\G_n$ is embedded, not trivial in $H_1(g_n^{-1}(\overline{\C }),\Z )$ and has period zero, Proposition~\ref{propos1} implies that $\G _n$ can be oriented so that
$[\G_n]=[\g_2(n)]$ in $H_1(g_n^{-1}(\overline{\C }),\Z )$.
Viewed in $H_1(M_n,\Z )$, the classes $[\G_n],[\g_2(n)]$ differ in a finite sum of loops around ends, from where
\begin{equation}
\label{eq:fluxGng2n}
F(\G _n)=i \overline{b_n}+t(n)\pi a_n,
\end{equation}
where $t(n)\in \Z $. Since both $\G_n$, $\g_2(n)$ can be chosen in the same fundamental domain of the doubly periodic lifting $\widetilde M_n$ of $M_n$ lying between two horizontal planes $\Pi ,\Pi +T_n$,
the embeddedness of both curves insures that $\{t(n)\}_n$ is bounded.
Passing to a subsequence, we can suppose that $t=t(n)$ does not depend on $n$.
Also note that $t$ is even since the periods of $M_n$ along both $\G_n,\g_2(n)$ are zero.
By item {\it 7} of Proposition~\ref{propos1}, $(F_{\G _n})_3=(F_{\g _2(n)})_3=2\pi $.
This property implies that the lengths of $L_1(n),L_2(n)$ diverge to $\infty$ as $n\to\infty$, see~\cite{mpr1} page~118. Those arguments in~\cite{mpr1} also show that the limit as $n\to\infty$ of the horizontal flux $F(\G_n)$ can be computed as the limit of the horizontal fluxes along $L_1(n)\cup L_2(n)$, and that the horizontal flux along $L_1(n)\cup L_2(n)$ divided by the length of $L(n)$ converges to a complex number of modulus 2 as $n\to\infty$, so $F(\G_n)\to\infty$.
Dividing equation (\ref{eq:fluxGng2n}) by $|F(\G _n)|$ and using that $b_n$ is bounded, we have (after extracting a subsequence) that both $|F(\G _n)|^{-1}F(\G _n)$, $t|F(\G _n)|^{-1}\pi a_n$ converge to the same limit $e^{i\t }$, $\t \in [0,2\pi )$, from where $t\neq 0$ and $a_n\to \infty $. Since $a_n\in\R^*$, we also have $\t =0$ or $\pi $. In particular, the direction of the segment $L(n)$ tends to the direction of the $x_2$-axis as $n\to \infty $.

Now consider respective liftings $\widehat L(n), \widehat{\cal D}_1(n), \widehat{\cal D}_2(n)$ of $L(n),{\cal D}_1(n),{\cal D}_2(n)$ such that $\widehat L(n)\cup\frac{1}{\l_n}\left[\widehat{\cal D}_1(n), \widehat{\cal D}_2(n)\right]$ lies in the same fundamental domain of $\widetilde M_n/T_n$.
Let $\widehat L^*(n)$ be the length minimizing horizontal segment from $\frac{1}{\l _n}\partial \widehat{\cal D}_2(n)$
to $\frac{1}{\l _n}\partial \widehat{\cal D}_1(n)+H_n$.
Let $\G^*_n\subset M_n$ be another embedded closed curve constructed in a similar way as $\G_n$, i.e. $\G^*_n=L^*_1(n)^{-1}*\be^*_1(n)*L^*_2(n)*\be^*_2(n)$ where $L^*_1(n),L^*_2(n)$ are liftings in consecutive sheets of $\widehat L^*(n)/H_n$ and $\be^*_1(n),\be^*_2(n)$ are small arcs inside $\frac{1}{\l _n}{\cal H}_1(n),\frac{1}{\l _n}{\cal H}_2(n)$ respectively.
As before, each of the $\be^*_i(n)$ consists of the consecutive union of three arcs, two at almost constant height joined by a central one which we choose as $\be_i(n)\cap N_n^{-1}(\esf^1)$, $i=1,2$.
We orient $\G^*_n$ in such a way that $\G _n,\G^*_n$ share their orientations along the arcs $\be_i(n)\cap N_n^{-1}(\esf^1)$, $i=1,2$.
Viewed in $H_1(g_n^{-1}(\overline{\C }),\Z )$, it holds $[\G _n]=-[\G^*_n]$, see Figure~\ref{magda2No} left.
As above, we have that after passing to a subsequence, $F(\G^*_n)=-i \overline{b_n}+t^*\pi a_n$ for certain nonzero even integer $t^*$.

\begin{figure}
\centerline{\includegraphics[width=15cm,height=6.16cm]{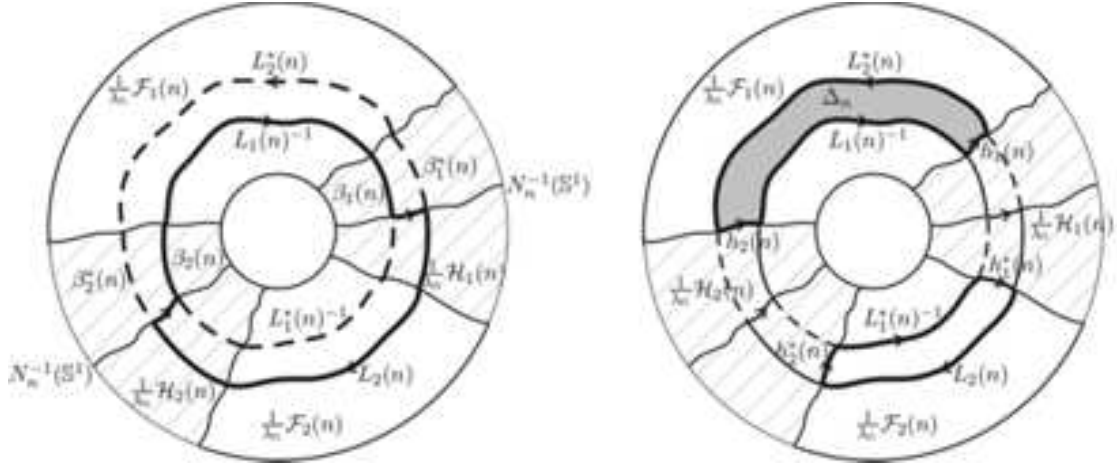}}
\caption{Left:
The curves $\G_n, \G_n^*$. Right: The curve $\zeta_n$.}
\label{magda2No}
\end{figure}

Let $\zeta _n\subset M_n$ be the embedded closed curve defined by $\zeta _n=L_1(n)^{-1}*h_1(n)*L^*_2(n)*h_2(n)$,
where $h_1(n)\subset \frac{1}{\l_n}\partial {\cal H}_1(n)$ is a helix-type curve joining the end point of $L_1(n)^{-1}$ with the starting point of $L^*_2(n)$ (so the covering map $f_n$ restricts to $h_1(n)$ as a diffeomorphism onto an arc in $\frac{1}{\l _n}\partial{\cal D}_1(n)$ arbitrarily close to a halfcircumference), and $h_2(n)\subset\frac{1}{\l_n}\partial {\cal H}_2(n)$ is similarly defined. Note that we can endow $\zeta_n$ with a natural orientation by restricting the orientation of either $\G_n$ or $\G_n^*$ to their common arcs with $\zeta_n$ (both restrictions give rise to the same orientation).
$\zeta_n$ is the boundary of a disk $\Delta_n$ inside $\frac{1}{\l_n}[{\cal F}_1(n)\cup{\cal F}_2(n)]$. We do not lose generality by supposing $\Delta_n\subset\frac{1}{\l_n}{\cal F}_1(n)$.
Since $f_n$ restricts to $\zeta_n$ as a diffeomorphism onto a curve that represents the period vector $H_n$, we deduce that $\Delta_n$ has just one end of $M_n$.
In particular $F(\zeta_n)=\ve\pi a_n$, with $\ve=\pm 1$.
Now consider for $i=1,2$ the helix-type curve $h^*_i(n)\subset \frac{1}{\l_n}\partial {\cal H}_i(n)$ that joins the corresponding extreme points of $L_1^*(n)^{-1}$ and of $L_2(n)$. $h_i^*(n)$ is almost opposite to $h_i(n)$ with respect
to the almost axis of $\frac{1}{\l_n}\partial {\cal H}_i(n)$, and let $\zeta^*_n$ be the closed curve $L^*_1(n)^{-1}*h^*_1(n)*L_2(n)*h^*_2(n)$ endowed with the orientation induced by the ones of either $\G_n$ or $\G_n^*$ along their common arcs with $\zeta^*_n$.
The same argument as before gives $F(\zeta^*_n)=\pm\pi a_n$.
Since the unit conormal vector along $L_1(n)^{-1}$ projects horizontally on the same side as the unit conormal vector along $L_2(n)$, we conclude that $F(\zeta^*_n)=F(\zeta_n)$.
A similar argument with the unit conormal vector implies that $\frac{F(\G_n)}{F(\G^*_n)}$ converges to a positive real number, from where $\frac{t}{t^*}>0$.

As we are assuming $\frac{r_n}{\l_n}\to 0$, then the lengths of $\be_i(n),\be^*_i(n),h_i(n),h^*_i(n)$ go to 0 as $n\to\infty$, $i=1,2$. Thus the limit of $F(\G _n)+F(\G^*_n)$ coincide with the limit of $F(\zeta_n)+F(\zeta^*_n)= 2\ve\pi a_n$ as $n\to\infty$, and
\[
2=\lim _{n\to \infty }\frac{\left( F(\G _n)+F(\G^*_n)\right) }{\ve \pi a_n}=
\lim _{n\to \infty }\frac{-2i\overline{b_n}+(t+t^*)\pi a_n}{\ve \pi a_n}=\ve (t+t^*),
\]
which contradicts that both $t, t^*$ are nonzero even integers with the same sign, thereby proving Proposition~\ref{proposcurvestim}.
\end{proof}

\begin{remark}
If we remove the hypothesis $(ii)$ in Proposition~\ref{proposcurvestim}, then it is possible to find two
highly sheeted vertical helicoids forming inside surfaces in ${\cal S}$, which makes their curvatures to blow-up.
For instance, the standard examples $M_{\t ,0,\pi /2}$ with $\t \nearrow \pi /2$ contain
two helicoids forming with axes joined horizontally by a line parallel to the period vector at the ends.
The surfaces $M_{\t ,0,\be }$ with $\be <\t <\pi/2$ and $\be \nearrow \pi /2$ have also two forming helicoids,
but now their axes join by a horizontal line orthogonal to the period vector at the ends.
\end{remark}

\section{Uniqueness of examples around the singly periodic Scherk surfaces.}
\label{secuniqaroundScherk1p}
In this Section we will prove that if $\{ M_n\} _n\subset \widetilde{\cal S}$ degenerates in a singly periodic
Scherk minimal surface (case {\it (iv)} of Proposition~\ref{proposdifferentlimits}), then $L(M_n)$
 tends to a tuple in $\C^{4k}$. In particular, the classifying map
$C:\widetilde{\cal S}\to \R^*\times \C$ cannot be proper. In order
to overcome this lack of properness, we will prove that only standard examples can occur in $\widetilde{\cal S}$ nearby
the singly periodic Scherk limit. This will be essential when proving that the restriction of $C$ to the space of
nonstandard examples is proper (Theorem~\ref{thproper}).

\begin{proposition}
\label{propos2kScherks1p}
 Let $\{M_n\}_n\subset \widetilde{\cal S}$ be a sequence of marked surfaces with $\{H_n\}_n\to
 H_\infty=(0,\pi a,0)$ ($a\in\R^*$), $\{T_n\}_n\to \infty$ (for any choice of $T_n$ as in Proposition~\ref{proposdifferentlimits}) and $\{C(M_n)\}_n \to
(a,b)$ where $b\in \C $. Then for $n$ large, the geometric surface $M_n$ is close to $2k$ translated images of
arbitrarily large compact regions of a singly periodic Scherk minimal surface of genus zero  with two
horizontal ends, together with $2k$ annular regions $C_n(1),\ldots ,C_n(2k)$ each of
which has two distinct simple branch points of the Gauss map $g_n:M_n\to \overline{\C }$.
 Moreover, there exists a nonhorizontal plane $\Pi\subset\R^3$ such that any
annulus $C_n(j)$ is a graph over the intersection of $\Pi/H_n$ with a certain
horizontal slab, $j=1,\ldots ,2k$.
\end{proposition}
\begin{proof}
Note that it suffices to show that there exists a subsequence of $\{ M_n\}_n$
verifying the conclusions of the Proposition.
Since the total branching number of the Gauss map $N_n$ of $M_n$ is fixed $4k$,
we can find a small $\varepsilon >0$ such that for every $n$, there exists a
$\theta_n\in \esf^2\cap \{ x_2=0\}$ with the disk $D_{\esf^2}(\theta_n,\varepsilon )$
of radius $\varepsilon$ in $\esf^2$ around $\theta_n$ disjoint from the branch locus of $N_n$ and from the North and South Poles.
In particular, $N_n^{-1}(\theta_n)$ consists of $2k$ distinct points $p_n(1),\ldots ,p_n(2k)\in M_n$. We label by $\widetilde{p}_n(i)$ the lift of $p_n(i)$ to a fundamental domain of the doubly periodic surface $\widetilde{M}_n
\subset \R^3$ whose quotient is $M_n$, $i=1,\ldots ,2k$.
Since $\{ C(M_n)\}_n$ converges in $\R^*\times \C $, Proposition~\ref{proposcurvestim} and Lemma~\ref{entornotub} imply that the sequence $\{ \widetilde{M}_n-\widetilde{p}_n(1)\}_n$ has uniform curvature and
local area bounds. After passing to a subsequence, $\widetilde{M}_n-\widetilde{p}_n(1)$ will converge
to a (not necessarily connected) properly embedded minimal surface $\widetilde{M}_{\infty }(1)\subset \R^3$ uniformly on
compact subsets of $\R^3$.
Since $\t_n\in\esf^2\cap\{x_2=0\}$, any plane orthogonal to the limit of $\{\t_n\}_n$ contains the vector $H_\infty$. Using this fact, minor modifications in the arguments inside the proof of Proposition~\ref{proposcurvestim} (when we proved that ${\cal H}_2$ is not flat) show that $\widetilde{M}_{\infty }(1)$ is not flat.
Applying Proposition~\ref{proposdifferentlimits} with the modified hypothesis of
the limit surface not being flat instead of the normalization of the absolute
curvature to have a maximum at the origin, and using that $H_n$ converges
while any choice of $T_n$ diverges,
we deduce that $\widetilde{M}_{\infty}(1)$ is a singly periodic Scherk minimal surface of genus zero, two of whose ends are horizontal.

Consider an infinite closed horizontal cylinder ${\cal C}\subset \R^3$ around the $x_2$-axis with radius large enough such that if we define $\Omega (1)=\big(\widetilde{M}_{\infty}(1)\cap {\cal C}\big)/H_{\infty }$, then $\big(\widetilde{M}_\infty(1)/H_{\infty}\big)-\Omega (1)$ consists of four extremely flat annular Scherk-type ends.
Since the region $\widetilde{\Omega}_n(1)=[\widetilde{M}_n-\widetilde{p}_n(1)]\cap {\cal C}$ satisfies that $\widetilde{\Omega}_n(1)/H_n$ converges uniformly to $\Omega(1)$ as $n\to \infty $ and $\Omega(1)$ has injective Gauss map, we deduce that  $\widetilde{\Omega}_n(1)/H_n$ has also injective Gauss map for $n$ large.
This implies that the points $\widetilde{p}_n(i)-\widetilde{p}_n(1)$, $2\leq i\leq 2k$, are outside $\widetilde{\Omega}_n(1)$.
Using the same arguments we conclude that each of the sequences $\{\widetilde{M}_n-\widetilde{p}_n(2)\}_n,\ldots ,\{\widetilde{M}_n-\widetilde{p}_n(2k)\}_n$ (after passing to a common subsequence) converges to respective singly periodic Scherk minimal surfaces $\widetilde{M}_{\infty }(2),\ldots ,\widetilde{M}_{\infty}(2k)$, each one with two horizontal ends. Note that these singly periodic Scherk limits have the same period vector $H_{\infty }$.
Taking the radius of ${\cal C}$ large enough, we can assume that $\big(\widetilde{M}_\infty(i)/H_{\infty}\big)-\Omega (i)$ consists of four extremely flat annular Scherk-type ends, where
$\Omega (i)=\big(\widetilde{M}_{\infty}(i)\cap {\cal C}\big)/H_{\infty }$ for each $i=1,\ldots ,2k$.
As before, if we define $\widetilde{\Omega}_n(i)=[\widetilde{M}_n-\widetilde{p}_n(i)]\cap {\cal C}$ then $\widetilde{\Omega}_n(i)/H_n$ converges uniformly to $\Omega (i)$ as $n\to \infty $. Furthermore, $\big[\big(\widetilde{\Omega}_n(1)+\widetilde{p}_n(1)\big)\cup \ldots \cup \big(\widetilde{\Omega}_n(2k)+\widetilde{p}_n(2k)\big)\big]/H_n$ embeds into $M_n$.
Since any compact horizontal level section of $M_n$ is connected and any compact horizontal level section of $\widetilde{\Omega}_n(i)/H_n$ is a closed curve,  we conclude that the minimum closed horizontal slab $S_n(i)\subset \R^3/H_n$ containing to $[\widetilde{\Omega}_n(i)+\widetilde{p}_n(i)]/H_n$ satisfies $S_n(i)\cap S_n(j)=\emptyset$ whenever $i\neq j$. Since the radius of the cylinder ${\cal C}$ can be made arbitrarily large, we also conclude that both the width of the slabs $S_n(i)$ and $(T_n)_3$ diverge to $\infty$ as $n$ increases.

We claim that all the limits $\widetilde{M}_{\infty }(1),\ldots ,\widetilde{M}_{\infty
}(2k)$ are in fact the same singly periodic Scherk minimal surface. Clearly it
suffices to check that the angle between the nonhorizontal ends of these surfaces and
the horizontal does not depend on $i=1\ldots ,2k$. Without loss of
generality, we can assume that $\widetilde{M}_{\infty }(1),\ldots
,\widetilde{M}_{\infty }(2k)$ are ordered increasingly in heights, in the sense that
for all $n$, $S_n(i+1)$ lies above $S_n(i)$, $i=1,\ldots ,2k-1$. Since the absolute
total curvature of $M_n$ is $8k\pi $, the one of a Scherk minimal surface with genus zero is $4\pi$ and we dispose of $2k$ of these limit surfaces, we conclude that the Gauss map $N_n$
restricted to the complement of $M_n\cap \left(\cup_{i=1}^{2k}(S_n(i)/T_n)\right)$
covers a set in $\esf^2$ of arbitrarily small area.
Let $C_n(j)$ be the component of $M_n-\left(\cup_{i=1}^{2k}(S_n(i)/T_n)\right)$
that glues to $\big[\widetilde\Omega_n(j)+\widetilde{p}_n(j)\big]/{\cal P}_n,\big[\widetilde\Omega_n(j+1)+\widetilde{p}_n(j+1)\big]/{\cal P}_n$.
Since $C_n(j)$ is a compact annulus in $M_n$, the Divergence Theorem shows that the flux vectors of $C_n(j)$ along its boundary curves are opposite; but such flux vectors converge as $n\to \infty$ to the fluxes of the limit Scherk surfaces $\widetilde{M}_{\infty }(j)/H_{\infty
},\widetilde{M}_{\infty }(j+1)/H_{\infty }$ around their respective upward and downward pointing nonhorizontal ends.
Thus, $\widetilde{M}_{\infty }(1)=\ldots =\widetilde{M}_{\infty }(2k)$ (and the Gauss map of $\widetilde{M}_{\infty }(j)$
is opposite to the one of $\widetilde{M}_{\infty }(j+1)$).

Since the Gauss map of a singly periodic Scherk surface of genus zero is unbranched, the $4k$
branch values (counting multiplicity) of the Gauss map $g_n$ of $M_n$ are located in
$C_n(1)\cup \ldots\cup C_n(2k)$. Given $j=1,\ldots,2k$, $g_n$ restricts to each of
the two boundary components of the annulus $C_n(j)$ as a bijection onto the boundary
of a small spherical disk centered at a point $N(j)$, this one being the limit normal vector of $\widetilde{M}_{\infty}(j)$ at its upward pointing nonhorizontal Scherk-type end or equivalently, the limit normal vector of $\widetilde{M}_{\infty }(j+1)$ at its downward pointing nonhorizontal Scherk-type end.
By gluing two suitable disks $D_1,D_2$ to $C_n(j)$ along its boundary
components, one can construct a meromorphic degree 2 map $G:C_n(j)\cup D_1\cup
D_2\to \overline{\C }$. Since $C_n(j)\cup D_1\cup D_2$ is a sphere, Riemann-Hurwitz
formula gives that $G$ has total branching number $2$ and so, $G$ has exactly
two distinct simple branch points which lie necessarily in $C_n(j)$. Finally, let $\Pi \subset \R^3$ be a
plane parallel to the nonhorizontal ends of $\widetilde{M}_{\infty }(1)$. Since for
$n$ large $g_n|_{C_n(j)}$ is contained in an arbitrarily small spherical disk
centered in $N(j)$, we conclude that $C_n(j)$ can be expressed as the graph of a
function $u_n(j):(\Pi/H_n)\cap S'_n(j)\to \R $, where $S'_n(j)$ is the horizontal
slab (quotiented by $H_n$) between $S_n(j)$ and $S_n(j+1)$.
This finishes the proof.
\end{proof}

Let $S_{\rho }$ be the singly periodic Scherk minimal surface that appears as a limit in
Proposition~\ref{propos2kScherks1p} (clearly we can assume $0<\rho \leq 1$). Since the period vector $H_{\infty }$
of $S_{\rho }$ points to the $x_2$-axis, the normal vectors at the ends of $S_{\rho }$ (stereographically projected) are $0,\infty, \rho ,-1/\rho$, so we can parametrize $S_{\rho }$
by the Weierstrass data
\begin{equation}
\label{eq:Scherk1pdata} g(z)=z,\quad dh=c\frac{dz}{(z-\rho )(\rho z+1)},\quad  z\in \overline{\C }-\{
0,\infty ,\rho ,-1/\rho \} ,
\end{equation}
where $c\in \R^*$. To determine $c$, note that the intersection of $S_{\rho }$ with the quotient by $H_{\infty }$ of
a horizontal plane at large positive height consists of a compact embedded curve $\G
$ with period $H_{\infty }$ (up to sign). $\G$ is the  uniform limit of compact horizontal level
sections $\G _n$ of the surfaces $M_n$ (suitably translated and with the notation of
Proposition~\ref{propos2kScherks1p}). Since the vertical parts of the period and flux vectors of
$M_n$ along $\G _n$ are respectively 0, $2\pi$  for all $n$, it follows that
\[
2\pi i=c\int_{\G }\frac{dz}{(z-\rho )(\rho z+1)}=2\pi ic \mbox{ Res}_\rho
\frac{dz}{(z-\rho )(\rho z+1)}=\frac{2\pi ic}{\rho^2+1},
\]
thus $c=\rho^2+1$ (we have assumed that the limit normal vector of $S_\rho$ at its upward pointing nonhorizontal end is $\rho$ and oriented $\G $ so that the second equality above holds; these choices determine the sign of $c$).

\subsection{Weierstrass data.}
Next we give a local chart for ${\cal W}$ around the boundary point described in
Proposition~\ref{propos2kScherks1p}.\label{Hurwitzschemes}
 Fix $\rho \in (0,1)$ and let $D(*,\ve )\subset \C$ be a small disk of radius $\ve >0$ centered at $*=\rho
,-1/\rho $. Given $k$ unordered couples of points $a_{2i-1},b_{2i-1}\in D(\rho ,\ve
)$ with $a_{2i-1}\neq b_{2i-1}$ and another $k$ couples $a_{2i},b_{2i}\in D(-1/\rho
,\ve )$ with $a_{2i}\neq b_{2i}$, $1\leq i\leq k$, we can construct a marked
meromorphic map $g\in {\cal W}$ associated to these couples as follows. Consider $2k$
copies $\overline{\C }_1,\overline{\C }_2,\ldots ,\overline{\C }_{2k}$ of
$\overline{\C }$. Cut $\overline{\C }_1$ along small disjoint arcs $\be_1,\be_2$ so that
$\be_1$ joins $a_1$ with $b_1$ and $\be_2$ joins $a_2$ with $b_2$. Cut $\overline{\C }_2$ along a copy of
$\be_2$ and also along an arc $\be_3$ joining $a_3$ with $b_3$, and glue
$\overline{\C }_1$ with $\overline{\C }_2$ along the common cut $\be_2$ in the usual way.
Repeat the process  so that $\overline{\C }_{2k}$ glues with
$\overline{\C }_1$ along the common cut $\be_1$. This surgery produces a torus $\M$ and the natural $z$-map on each copy of $\overline{\C }$
gives a well defined degree $2k$ meromorphic map $g:\M \to \overline{\C }$ with
branch values $\{a_1,\ldots,b_{2k}\} $. The ordered list of zeros and poles of
the marked meromorphic map to be defined will be $(0_1,\ldots
,0_{2k},\infty_1,\ldots ,\infty_{2k})$, where the subindexes refer to the copy of $\overline{\C }$
which the zero or pole of $z$ belongs to (we do not lose generality by assuming that $0,\infty$ do not lie in none of the $\be_j$-curves). Finally, the nontrivial homology class
$[\g ]\in H_1(\M -\{ \mbox{zeros, poles}\} ,\Z )$ is defined to be the class of $\{ |z|=1\}$ in $\overline{\C }_1$
with the anticlockwise orientation.

\begin{remark}
\label{remarkcurveScherk1p}
For any $\rho \in (0,1)$, the circle $\{|z|=1\}$ is a closed embedded curve in $S_{\rho}$ with period zero. When
$\rho =1$ we must perturb slightly $\{ |z|=1\} $ to keep the vanishing period condition true.
This fact justifies the above choice of the homology class $[\g ]$.
\end{remark}

With the above procedure, the map $(a_1,b_1,\ldots ,a_{2k},b_{2k})\mapsto g\in {\cal W}$ is not
injective, as one can exchange $a_i$ by $b_i$ obtaining the same $g$. Rather than parametrizing these marked meromorphic maps by the lists of its branch values, one can use the symmetric elementary polynomials in two
variables for each couple $a_i,b_i$. Note that when $a_i,b_i$ are close to $\rho $
(resp. $-1/\rho $) then $(a_i+b_i,a_ib_i)$ is close to $(2\rho ,\rho^2)$ (resp.
$(-2/\rho ,1/\rho^2)$). The computations that follow simplify if we use the
arithmetic and geometric means instead of the elementary
symmetric polynomials, which is possible since the map $(u,v)\mapsto (u/2,\sqrt{v})$ is a local
diffeomorphism outside $(0,0)$. We define
\[
x_i=\frac{1}{2}(a_i+b_i),\qquad y_i=\sqrt{a_ib_i},
\]
so $(x_i,y_i)$ lies in a neighborhood of $(\rho ,\rho )$ or $(-1/\rho ,1/\rho )$. Given $\ve
'>0$, we label ${\cal U}(\ve')=\left[D(\rho ,\ve ')\times D(\rho ,\ve ' )\times
D({\textstyle \frac{-1}{\rho }},\ve ')\times D({\textstyle \frac{1}{\rho }},\ve
')\right]^k$. As $x_i^2-y_i^2=\frac{1}{4}(a_i-b_i)^2$, the condition $a_i\neq b_i$ necessary for the above
construction is equivalent to $x_i^2\neq y_i^2$. Let ${\cal A}=\{(x_1,y_1,\ldots ,x_{2k},y_{2k})\ $ $| \ x_i^2=y_i^2\ \mbox{for some }i=1\ldots
,2k\} $. Clearly ${\cal A}$ is an analytic subvariety of $\C^{4k}$.  It can be shown that for $\ve '>0$ small, the correspondence ${\bf z}=(x_1,y_1,\ldots ,x_{2k},y_{2k})\in {\cal U}(\ve ')-{\cal A}\mapsto \aleph ({\bf z})=g\in{\cal W}$ defines a local chart for ${\cal W}$.

\begin{remark}
\begin{description}
\item[]
\item[{\it (i)}]\label{remarkordering}
If a marked meromorphic map $g=\aleph ({\bf z})$ produces a marked surface $M$, then the ordered list
$(0_1,\ldots ,0_{2k},\infty_1,\ldots ,\infty_{2k})$ does not necessarily coincide with the ordering on the ends of $M\in \widetilde{\cal S}$. This is only a matter of notation and will not affect to the arguments that follow.
\item[{\it (ii)}]\label{notalimscherk1p}
Consider a sequence $\{ M_n\}_n\subset \widetilde{\cal S}$ with $\{C(M_n)=(a_n,b_n)\}_n\to(a,b)\in\R^*\times\C$ and $\{T_n\}_n\to\infty$ as $n\to\infty$ (for any choice of $T_n$).
By Proposition~\ref{propos2kScherks1p}, the sequence of geometric surfaces $\{M_n\}_n$ converges uniformly to $2k$-copies of a singly periodic Scherk minimal surface $S_\rho$ parametrized as in (\ref{eq:Scherk1pdata}), for certain $\rho \in(0,1]$.
Let $\Gamma$ be the curve $\{|z|=1\}$ viewed in one of the copies of $S_\rho$.
Clearly, $\Gamma$ is the uniform limit as $n\to\infty$ of a sequence of closed curves $\Gamma_n\subset M_n$ with $P_{\Gamma_n}=0$.
After exchanging the homology class of the marked surface $M_n$ by $[\G_n]\in H_1(M_n,\Z)$, we can see the same geometric surface $M_n$ as a new marked surface $M'_n$ inside the domain of the chart $\aleph$ for $n$ large enough (also note that the second component of $C(M'_n)$ differs from $b_n$ in a fixed even multiple of $i \pi a_n$).
\end{description}
\end{remark}

\subsection{Holomorphic extension.}

When ${\bf z}\in {\cal A}$, the continuous extension of the above cut-and-paste process gives a
Riemann surface with nodes (see~\cite{Ima} page 245 for the definition of a Riemann surface with nodes),
each node occurring between copies $\overline{\C }_{j-1},\overline{\C }_j$ such that $a_j=b_j$.
The corresponding differential $\phi $ also extends through ${\bf z}$, in the following manner.
\begin{proposition}
\label{proposholomextensionscherk1p}
Each ${\bf z}\in {\cal A}$ produces $l$ spheres $S_1,\ldots ,S_l$ joined
by $l$ node points $P_i,Q_i\in S_i$ (here $Q_i=P_{i+1}$ and the subindexes are cyclic),
$g$ degenerates in $l$ nonconstant meromorphic maps $g(i):S_i\to \overline{\C }$ with
the degrees of $g(1),\ldots ,g(l)$ adding up to $2k$ and $g(i)$ takes the values $\rho$ and/or $-1/\rho$ at $P_i,Q_i$.
$\phi $ degenerates in the $l$ unique meromorphic
differentials $\phi (i)$ on $S_i$, such that $\phi (i)$ has exactly two simple poles at $P_i,Q_i$ with residues $1$
at $P_i$ and $-1$ at $Q_i$ (these residues are determined by the equation $\int _{|z|=1}\phi =2\pi i $).
Finally, both $g$ and $\phi $ depend holomorphically on all parameters (including at points of ${\cal A}$).
\end{proposition}
\begin{proof}
See Lemma~8 of~\cite{mpr1} for a similar situation as ours (tori degenerating in spheres); see also
Section 3.4 of~\cite{tra1} for arbitrary genus.
\end{proof}

\begin{lemma}
\label{lemaLextendsScherk1p}
The ligature map $L$ extends holomorphically to ${\cal U}(\ve ')$.
\end{lemma}
\begin{proof}
Since ${\cal U}(\ve ')\cap {\cal A}$ is an analytic subvariety of ${\cal U}(\ve ')$, it
suffices to check that $L$ is bounded by the Riemann extension
Theorem for several variables (see for instance~\cite{GriHar} page~9).
Let $g_n=\aleph ({\bf z}(n))$ be a sequence with ${\bf z}(n)\in {\cal U}(\ve ')-{\cal A}$ converging to
${\bf z}\in{\cal U}(\ve ')\cap {\cal A}$.
Note that each component of $L(g_n)$ can be written as an integral along a curve $\a $ independent of $n$ which lies in one of the twice punctured spheres $S_j-\{ P_j,Q_j\} $ of Proposition~\ref{proposholomextensionscherk1p}, of a holomorphic differential $\varphi (n)$, and that $\{ \varphi (n)\} _n$ converge uniformly on $\a $ as $n\to \infty $ to some holomorphic differential $\varphi $ on $S_j-\{ P_j,Q_j\} $. From here we directly deduce that $L(g_n)$ is bounded as desired.
\end{proof}

\subsection{Partial derivatives and Inverse Function Theorem.}

\begin{theorem}
\label{teorLbiholScherk1p}
 There exists $\ve '>0$ small such that $L|_{{\cal U}(\ve
')}$ is a biholomorphism.
\end{theorem}
\begin{proof}
Given a list ${\bf z}=(x_1,y_1,\ldots,x_{2k},y_{2k})\in {\cal U}(\ve')-{\cal A}$, we will denote by
$g=\aleph ({\bf z})$ the associated  marked meromorphic map. Recall that $0_j,\infty _j$ belong to the copy
$\overline{\C }_j$ of $\overline{\C }$. For $j=1,\ldots,k$, we let $\G_{2j-1}$ denote the closed curve in
$g^{-1}(\overline{\C })$ that corresponds to the loop $\{|z|=1\}$ on $\overline{\C }_{2j-1}$ (if $\rho =1$, then we modify slightly $\{ |z|=1\} $ as mentioned in Remark~\ref{remarkcurveScherk1p}).
We orient $\G _1$ to coincide with the last component of $g$ and the reamining $\G_{2j-1}$ to be homologous to
$\G _1$ in $g^{-1}(\overline{\C })$. Hence
Res$_{0_{2j}}(\frac{\phi}{g})+\mbox{Res}_{0_{2j+1}}(\frac{\phi}{g})=\frac{1}{2\pi
i}\big(\int_{\G_{2j+1}}\frac{\phi}{g}-\int_{\G_{2j-1}}\frac{\phi}{g}\big)$ and a similar formula holds for poles
of $g$. Thus, the composition of $L$ with a certain regular linear
transformation in $\C^{4k}$ writes as $\widehat L:{\cal U}(\ve')\to\C^{4k}$ where
\[
\widehat L({\bf z})=\Big( \underbrace{\mbox{Res}_{0_{2j-1}}\Big( \frac{\phi }{g}\Big) ,
\mbox{Res}_{\infty _{2j-1}}(g\phi)}_{1\leq j\leq k},
\underbrace{\int_{\G_{2j-1}}\frac{\phi }{g}, \int_{\G_{2j-1}}g\phi}_{1\leq j\leq
k}\Big).
\]
By the Inverse Function Theorem, it suffices to prove that the differential of
$\widehat L$ at the point $(\rho ,\rho,\frac{-1}{\rho
},\frac{1}{\rho })^{k}\in \C^{4k}$ corresponding to ($2k$ copies of) the singly periodic Scherk minimal
surface $S_{\rho }$, is an automorphism of $\C^{4k}$. The
computations that follow are similar to those in the proof of Lemma 9 in
\cite{mpr1}, so we will only explain them briefly.

Fix $j=1,\ldots,k$. To compute $\frac{\partial \widehat L}{\partial x_{2j-1}}(S_\rho)$,
we differentiate at $x=\rho$ the composition of $\widehat L$ with the curve $x\in D(\rho,\ve')\mapsto
(\rho,\rho,\frac{-1}{\rho},\frac{1}{\rho},\ldots, \frac{1}{\rho},x,\rho,\frac{-1}{\rho},
\frac{1}{\rho},\ldots, \rho,\rho,\frac{-1}{\rho},\frac{1}{\rho})\in {\cal U}(\ve ')$, where $x$ is placed at
the $(4j-3)$-th component of the $4k$-tuple.
This curve produces $2k-1$ spheres $S_1, S_2,\dots,S_{2j-2},S_{2j},\ldots,S_{2k}$ and
meromorphic maps $g_1, g_2,\dots,g_{2j-2},g_{2j},\ldots,g_{2k}$ on them such that
\begin{enumerate}
\item For each $m=1,\dots,2k$ with $m\neq 2j-2,2j-1$, $g_m:S_m\to\overline{\C }$ is a biholomorphism.
Thus we can parametrize $S_m$ by $\overline{\C }$ with $g_m(z)=z$. With this
parameter, the node points correspond to $\rho$, $\frac{-1}{\rho}$, and $\phi=c\frac{dz}{(z-\rho)(z+\rho^{-1})}$ where $c\in\C^*$ is determined by the equation
$\int_{|z|=1}\phi=2\pi i$ (in particular, $c$ does not depend on $x$).
\item $g_{2j-2}:S_{2j-2}\to\overline{\C }$ has degree two. We can parametrize $S_{2j-2}$ by $\{(z,w)
\in \overline{\C }^2\ |\ w^2=z^2-2xz+\rho^2\}$, hence $w=\sqrt{z^2-2xz+\rho^2}$ is
well defined on $S_{2j-2}$ (we fix the sign of the square root so that $w\sim z-\rho $ in $\overline{\C}_{2j-1}$
and $w\sim -(z-\rho )$ in $\overline{\C}_{2j-2}$). Let $Q$ denote the point $z=-1/\rho $ in $\overline{\C }_{2j-1}$.
Then $g(z,w)=z$ and $\phi = c(x)\frac{dz}{(z+\rho^{-1})w}$, where
\begin{equation}
\label{c(x)2j-1Scherk}
-1=\Res _Q\phi =\Res _Q\frac{c(x)\, dz}{(z+\rho^{-1})w}=\frac{c(x)}{w(Q)}.
\end{equation}
\end{enumerate}
\par
\noindent
Since the components of $\widehat L$ computed on $S_m$ with $m\neq 2j-2,2j-1$ do
not depend on $x$, the corresponding derivative with respect to $x$ vanishes.
From (\ref{c(x)2j-1Scherk}) one has $c(\rho)=-w(Q)|_{x=\rho }=\frac{\rho^2+1}{\rho}$,
$c'(\rho)=-\left. \frac{d}{dx}\right| _{x=\rho }w(Q)=\frac{1}{\rho^2+1}$. An elementary calculation gives
\[
\begin{array}{ll}
{\displaystyle \left. \frac{d}{dx}\right|_{x=\rho }\mbox{Res}_{0_{2j-1}}\Big(\frac{\phi}{g}\Big) =
\left. \frac{d}{dx}\right|_{\rho }\mbox{Res}_{(0,-\rho )}\frac{c(x)\, dz}{z(z+\rho^{-1})w}=-c'(\rho ),}
&
{\displaystyle \left. \frac{d}{dx}\right|_{\rho }\mbox{Res}_{\infty _{2j-1}}(g\phi )=-c'(\rho ),}
\\
{\displaystyle \left. \frac{d}{dx}\right|_{\rho }\int_{\G_{2j-1}}\frac{\phi }{g}=\left. \frac{d}{dx}\right|_{\rho }\int_{|z|=1}\frac{c(x)\, dz}{z(z+\rho^{-1})w}=0,}
&
{\displaystyle \left. \frac{d}{dx}\right|_{\rho }\int_{\G_{2j-1}}g\phi =\frac{2\pi i }{\rho^2+1},
}
\end{array}
\]
where $\G_{2j-1}\subset S_{2j-2}$ is the connected lift of $\{|z|=1\}$ to
$\overline{\C }_{2j-1}$ through the $z$-map. Thus,
\[\begin{array}{rcccccccc}
\frac{\partial \widehat L}{\partial x_{2j-1}} = & \frac{1}{\rho^2+1} &
\Big(0,\ldots,0\hspace{-.3cm} &
,\hspace{-.1cm} \stackrel{(2j-1)}{-1},&
\hspace{-.3cm}\stackrel{(2j)}{-1},&
\hspace{-.3cm}0\, ,\ldots,0\hspace{-.3cm}&
\stackrel{(2k+2j-1)}{,\hspace{.5cm}0\hspace{.5cm,}}&
\hspace{-.3cm}\stackrel{(2k+2j)}{2\pi i},&
\hspace{-.2cm}0\,,\ldots,0\Big) .
\end{array}
\]
Analogous derivation of the composition of $\widehat{L}$ with similar curves gives
\[\begin{array}{rcccccccc}
 & & &_{(2j-1)} &\hspace{-.2cm}_{(2j)} & &\hspace{-.5cm}_{(2k+2j-1)} &\hspace{-.3cm}_{(2k+2j)} & \\

\frac{\partial \widehat L}{\partial y_{2j-1}} = \hspace{-.3cm}&
\frac{\rho^2}{\rho^2+1} \hspace{-.4cm}& \Big(0,\dots,0, \hspace{-.3cm}&
1-\frac{(\rho^2+1)^2}{\rho^4},&
\hspace{.2cm}-1\hspace{.5cm},&
\hspace{-.3cm}0,\dots,0 , &
\hspace{.2cm} 0\hspace{.5cm},&\hspace{-.3cm} 2\pi i&
\hspace{-.4cm},0,\dots,0\Big);\\

\rule{0cm}{.75cm}
\frac{\partial \widehat L}{\partial x_{2j}} = \hspace{-.3cm}&
\frac{\rho^2}{\rho^2+1} \hspace{-.4cm}& \Big(0,\ldots,0 ,\hspace{-.5cm}&
\hspace{.7cm}1\hspace{.8cm}, & \hspace{.5cm}1\hspace{.4cm}, &
\hspace{-.3cm}0,\ldots,0 , &
2\pi i\hspace{.4cm}, & 0\hspace{.4cm} &
\hspace{-.4cm},0,\ldots,0\Big);\\

\rule{0cm}{.75cm}\frac{\partial \widehat L}{\partial y_{2j}} = \hspace{-.3cm}&
\frac{(\rho^2+1)^2-1}{\rho^2+1} \hspace{-.4cm}& \Big(0,\dots,0, \hspace{-.3cm}&
\hspace{.7cm}1\hspace{.8cm}, & \hspace{-.3cm}\frac{1}{1-(\rho^2+1)^2}, &
\hspace{-.3cm}0,\dots,0 , &
2\pi i\hspace{.4cm}, & 0\hspace{.4cm} &
\hspace{-.4cm},0,\dots,0\Big).\\
\end{array}
\]
The absolute Jacobian of $\widehat L$ at the point $S_\rho$ is the absolute value of
the  determinant of the above files, which turns out to be $(2\pi)^{2k}$, and the
Theorem follows.
\end{proof}

\section{Uniqueness of examples around the catenoid.}
\label{secuniqaroundcat}
When a sequence $\{ M_n\}_n\in \widetilde{\cal S}$ degenerates in a
vertical catenoid (case {\it (i)} of Proposition~\ref{proposdifferentlimits}), the
residues in the ligature map $L$ diverge to $\infty $. In this Section we will
modify $L$ to have a well defined locally invertible extension through this boundary point of
${\cal W}$.

\begin{proposition}
\label{propos2kcatenoids}
 Let $\{M_n\}_n\subset \widetilde{\cal S}$ be a sequence with $\{C(M_n)=(a_n,b_n)\}_n \to
(\infty ,0)$. Then for $n$ large, the geometric surface $M_n$ is close to $2k$ translated images of arbitrarily
large compact regions of a catenoid with flux $(0,0,2\pi )$, together
with $2k$ regions $C_n(1),\ldots ,C_n(2k)$. Each $C_n(j)$ is a twice punctured annulus with one left end,
one right end of $M_n$ and two distinct simple branch points of its Gauss map. Furthermore, $C_n(j)$ is a
graph over its horizontal projection on $\{ x_3=0\} /H_n$.
\end{proposition}
\begin{proof}
As in Proposition~\ref{propos2kScherks1p}, it suffices to check that there exists a
subsequence of $\{ M_n\}_n$ verifying our assertions.
Given $n\in \N$, let $p_n(1),\ldots ,p_n(2k)\in M_n$ be the $2k$ distinct points applied
by the Gauss map $N_n$ of $M_n$ on a prescribed regular value $\t_n\in
\esf^2\cap \{ x_2=0\}$ such that a spherical disk $D_{\esf^2}(\t_n,\ve )$ of sufficiently small
uniform radius $\ve$ centered at $\t_n$ does not contain neither branch values of $N_n$ nor the
North and South Poles. Let $\widetilde{M}_n\subset \R^3$
be the doubly periodic surface obtained by lifting $M_n$ and $\widetilde{p}_n(i)$ the point that corresponds
to $p_n(i)$ in a fixed fundamental domain of $\widetilde{M}_n$, $i=1\ldots ,2k$. By
Proposition~\ref{proposcurvestim} and Lemma~\ref{entornotub}, the sequence $\{
\widetilde{M}_n-\widetilde{p}_n(1)\}_n$ has uniform curvature and local area bounds.
After extracting a subsequence, $\widetilde{M}_n-\widetilde{p}_n(1)$ converges
uniformly on compact subsets of $\R^3$ to a properly embedded minimal surface
$\widetilde{M}_{\infty }(1)\subset \R^3$. The same argument as in the proof of
Proposition~\ref{propos2kScherks1p} proves that $\widetilde{M}_{\infty }(1)$ is not flat,
so it lies in one of the six cases in Proposition~\ref{proposdifferentlimits}. Using that
$a_n\to\infty$ as $n\to \infty $, we discard the cases {\it (iv), (v), (vi)} of that Proposition.
Since the vertical part of the flux of all the $M_n$ along a compact horizontal section is $2\pi$ and the vertical part of the conormal vector to $M_n$ along such a horizontal section cannot vanish, we conclude that $\widetilde{M}_{\infty }(1)$ is not a vertical helicoid.
Now suppose that $\widetilde{M}_{\infty }(1)$ is a Riemann minimal example, and let $F\in\C^*$ be the horizontal part of its flux along a compact horizontal section. Then for $n$ large there exists a closed curve
$\G_n\subset\widetilde{M}_n-\widetilde{p}_n(1)$ such that $F(\G_n)$ tends to $F$ as $n\to \infty $.
But $F(\G_n)-i\overline{b_n}=t(n)\pi a_n$ for certain even integer $t(n)$, which gives a contradiction after taking limits.
Hence $\widetilde{M}_{\infty }(1)$ is a vertical catenoid with flux $(0,0,2\pi)$.

Reasoning as in the proof of Proposition~\ref{propos2kScherks1p} we conclude that
for $n$ large, the points $\widetilde{p}_n(i)-\widetilde{p}_n(1)$, $2\leq i\leq 2k$, are outside a
large compact domain of $\widetilde{M}_n-\widetilde{p}_n(1)$ arbitrarily close to a
vertical catenoid, and that after passing to a subsequence, $\{
\widetilde{M}_n-\widetilde{p}_n(2)\}_n,\ldots ,\{
\widetilde{M}_n-\widetilde{p}_n(2k)\}_n$ converge to different translations of
$\widetilde{M}_{\infty }(1)$. Since an arbitrarily large compact region $\Omega (i)$ of $\widetilde{M}_{\infty }(1)$ can be uniformly approximated by compact regions $\widetilde{\Omega}_n(i)\subset \widetilde{M}_n-\widetilde{p}_n(i)$, it follows that $(T_n)_3\to \infty$ as $n\to \infty $. Also note that the regions $[\widetilde{\Omega}_n(i)+\widetilde{p}_n(i)]/H_n$ can be chosen as the intersection of $\widetilde{M}_n/H_n$ with disjoint horizontal slabs $S_n(i)\subset\R^3/H_n$ whose widths go to $\infty$ as $n\to \infty $.

For $j=1,\ldots ,2k$, let $C_n(j)$ be the component of $M_n-\left(\cup_{i=1}^{2k}(S_n(i)/T_n)\right)$ that glues to $[\widetilde{\Omega}_n(j)+\widetilde{p}_n(j)]/{\cal P}_n$, $[\widetilde{\Omega}_n(j+1)+\widetilde{p}_n(j+1)]/{\cal P}_n$ (we can assume that $S_n(j+1)$ is directly above $S_n(j)\ $).
A straightforward modification of the argument in Proposition~\ref{propos2kScherks1p} using the
injectivity of the Gauss map of a catenoid shows that each $C_n(j)$ can be compactified by adding two ends
$e_1(n),e_2(n)$ of $M_n$ to obtain a compact annulus with two single branch points of the Gauss map $N_n$ of $M_n$, and that $N_n$ applies bijectively each boundary curve of $C_n(j)$ on the boundary of a small
spherical disk centered at $(0,0,\pm 1)$. As the fluxes around
$e_1(n),e_2(n)$ cancel (by the Divergence Theorem and because the flux of $M_n$ along each boundary curve of $C_n(j)$ tends to finite vertical but the flux at a Scherk-type end of $M_n$ is horizontal and arbitrarily large),
we deduce that $e_1(n)$ is a right end of $M_n$ and $e_2(n)$ a left end (or vice versa), both with the same limit normal vector. Finally, the horizontal projection of $C_n(j)$ onto its image in $\{ x_3=0\}
/H_n$ extends smoothly across $e_1(n),e_2(n)$ giving rise to a proper local diffeomorphism hence a finite covering map,
whose degree is one because $C_n(j)$ has one left end and one right end. This gives the graphical property of
$C_n(j)$ and finishes the proof of the Proposition.
\end{proof}

\subsection{Weierstrass data.}
Following the line of arguments in Section~\ref{secuniqaroundScherk1p}, we next show
a local chart for ${\cal W}$ around the boundary point that appears in
Proposition~\ref{propos2kcatenoids}. Given $i=1,\ldots ,k$, choose points
$a_{2i-1},b_{2i-1}$ (resp. $a_{2i},b_{2i}$) in a small punctured neighborhood of $0$ (resp.
of $\infty $) in $\overline{\C }$, such that $a_j\neq b_j$ for any $j$. These unordered couples can be considered as the branch values of a meromorphic map $g$ of degree $2k$, and a cut-and-paste construction analogous to the one in Section~\ref{secuniqaroundScherk1p} gives rise to a
marked meromorphic map $g\in {\cal W}$.
Since the roles of $a_j$ and $b_j$ are symmetric,
their elementary symmetric functions are right parameters in this setting. We introduce the parameters
\[
\begin{array}{lll}
x_j=\frac{1}{2}(a_j+b_j),&  y_j=a_j b_j &\mbox{if $j$ is odd}\\
x_j=\frac{1}{2}\left(\frac{1}{a_j}+\frac{1}{b_j}\right),&
y_j=\frac{1}{a_j b_j}& \mbox{if $j$ is even}
\end{array}
\]
so all parameters $x_j,y_j$, $1\leq j\leq 2k$, are close to $0$.
Also the conditions on $a_j$, $b_j$ translate into $y_j\neq x_j^2$ and $y_j\neq 0$.
In what follows,
we abbreviate ${\bf x}=(x_1,\ldots ,x_{2k}),\ {\bf y}=(y_1,\ldots ,y_{2k})$. Given
$\ve >0$, we let
\[
\begin{array}{l}
D(0 ,\ve )^{4k}=\{ ({\bf x},{\bf y})\in \C^{4k}\ | \ |x_j|,|y_j|<\ve \mbox{ for all }j=1\ldots ,2k\} ,\\
{\cal B}=\{ ({\bf x},{\bf y})\in D(0,\ve )^{4k}\ | \ x_j^2=y_j \mbox{ for some }j\},\\
\widehat{\cal B}=\{ ({\bf x},{\bf y})\in D(0,\ve )^{4k}\ | \ y_j=0 \mbox{ for some }j\} .
\end{array}
\]
${\cal B}\cup \widehat{\cal B}$ is an
analytic subvariety of the polydisk $D(0,\ve )^{4k}$ and the map $({\bf x},{\bf
y})\in D(0,\ve )^{4k}-({\cal B}\cup \widehat{\cal B})\mapsto \chi({\bf x},{\bf
y})= g\in{\cal W}$ is a local chart for ${\cal W}$. We also let $({\bf 0},{\bf 0})$ be the tuple
$({\bf x},{\bf y})=(0,\ldots ,0)\in D(0,\ve )^{4k}$.

\begin{remark}\label{notalimcat}
An argument as in Remark~\ref{notalimscherk1p}-{\it (ii)} shows that given a sequence $\{ M_n\}_n\subset \widetilde{\cal S}$ with $C(M_n)\to (\infty ,0)$, there exists another sequence of marked surfaces $\{M'_n\}_n$ inside the image of the chart $\chi$ such that for each $n$, $M_n,M'_n$ only differ in the homology class in the last component of the marked surface.
\end{remark}

\begin{figure}
\centerline{\includegraphics[width=11.1cm,height=3.5cm]{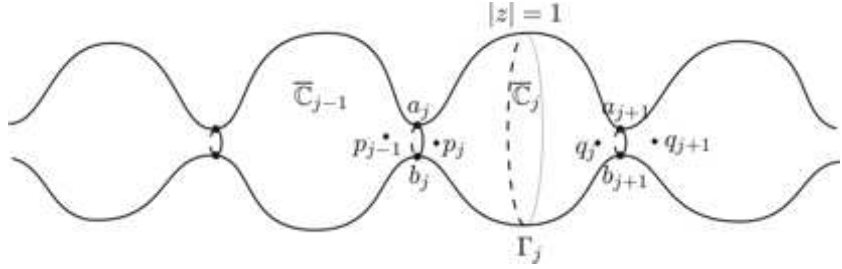}}
\caption{In
this case, $j$ is odd.}
\label{magda3}
\end{figure}

Next we find the equations to solve in order to produce an immersed minimal surface.
Let $\overline{\C}_j$ be the $j^{\mbox{th}}$ copy of $\overline{\C}$ and
let $\Gamma_j$ be the circle defined by $|z|=1$ in $\overline{\C}_j$, with
the positive orientation if $j$ is odd and the negative orientation
if $j$ is even. All these curves are homologous in $g^{-1}(\overline{\C })$, and $[\G _1]$ is the
last component of the marked meromorphic map $g$.
We also write $0_j$, $\infty_j$ for the points $z=0$ and $z=\infty$ in $\overline{\C}_j$,
see Figure~\ref{magda3} for the case $j$ odd.
Recall that $\phi$ is defined as the unique holomorphic 1-form with
$\int_{\gamma}\phi=2\pi i$. Define for $1\leq j\leq 2k$
\[
A_j=\left\{
\begin{array}{ll}
\int_{\Gamma_j}g^{-1}\phi &\mbox{($j$ odd)}\\
\int_{\Gamma_{j+1}}g\phi &\mbox{($j$ even)}
\end{array}
\right.
\qquad
B_j=\left\{
\begin{array}{ll}
\Res_{0_{j-1}} \left(g^{-1}\phi \right) \cdot \Res_{0_j} \left( g^{-1}\phi \right) &\mbox{($j$ odd)}\\
\Res_{\infty_{j-1}}\left( g\phi\right) \cdot \Res_{\infty_j}\left( g\phi \right) &\mbox{($j$ even)}
\end{array}
\right.
\]
In this definition and in the sequel we will adopt a cyclic convention on the subindexes, so when $j=1$, $j-1$ must be understood as $2k$. By the residue theorem, we have for odd $j$
\[
A_j-A_{j-2}=2\pi i \left(\Res_{0_j}g^{-1}\phi+\Res_{0_{j-1}}g^{-1}\phi\right) ,
\]
\[
A_{j+1}-A_{j-1}=2\pi i\left(\Res_{\infty_{j+1}}g\phi+\Res_{\infty_{j}}g\phi\right) .
\]
It follows easily that $g$ closes periods if and only if there exist $a\in\R$, $b\in\C$ such that
\begin{equation}
\label{eq:gclosesperiods}
\left.
\begin{array}{ll}
A_{2i-1}=b, A_{2i}=\overline{b} &\mbox{for all }i=1,\ldots ,k ,\\
B_j=-a^2 & \mbox{for all }j=1,\ldots ,2k,
\end{array}
\right\}
\end{equation}
which are the equations we have to solve.

\begin{remark}
\label{remarkResmultivalued}
The definition of $A_j$ and $B_j$ are motivated by the following fact: Fix $j$, say odd. There is no natural way
to distinguish between the two zeros of $g$ which are close to the branch points $a_j$, $b_j$.
We call them $0_{j-1}$ and $0_j$ for convenience, but this is a little bit artificial. Indeed the fact that one
zero belongs to $\overline{\C}_j$ or $\overline{\C}_{j-1}$ depends on the choice of the common cut $\beta_j$
connecting the points $a_j$ and $b_j$, but there is no way to choose $\beta_j$ depending continuously on the
parameters (this is a homotopy issue). In other words, the residue of $g^{-1}\phi$ at $0_j$ is not
a well defined function of the parameters. On the other hand, the unordered pair $\{0_{j-1},0_j\}$ depends
continuously on the parameters. For this reason, we should consider the elementary symmetric functions
of the residues at $0_{j-1}$ and $0_j$. This is essentially what we do in the definition of
$A_j$ and $B_j$.
\end{remark}

\subsection{Holomorphic extension.}
When $a_j=b_j$ (which corresponds to $y_j=x_j^2$) for a given $j$,
the definition of $g$ gives a Riemann surface with a node between $\overline{\C }_{j-1}$ and $\overline{\C }_j$.
More precisely, each $(\x ,\y )\in {\cal B}$ gives rise to a Riemann surface with nodes which consists of $l$ spheres
$S_j$ joined by node points $P_j,Q_j$ so that $P_j=Q_{j+1}$, $l$ nonconstant meromorphic maps $g(j):S_j\to
\overline{\C }$ with $\sum _j\deg (g(j))=2k$ and $g(j)(\{ P_j,Q_j\} )\subset \{ 0,\infty \} $, and
$l$ meromorphic differentials $\phi (j)$ on $S_j$ with just two simple poles at $P_j,Q_j$ and residues $1$ at $P_j$,
$-1$ at $Q_j$.

If $a_j=0$ and $b_j\neq 0$ (or vice versa) for $j$ odd, then the conformal structure between the copies
$\overline{\C }_{j-1}$ and $\overline{\C }_j$ does not degenerate, but the corresponding $z$-map has a double zero.
For $j$ even we have a similar behavior exchanging zero by pole.
Thus each $(\x ,\y )\in \widehat{\cal B}-{\cal B}$ produces a conformal torus $\M $ and a single meromorphic
degree $2k$ map $g:\M \to \overline{\C }$ with at least a double zero or pole. $\phi $ extends to the (unique) holomorphic differential on $\M$ with $\int _{\G _1}\phi =2\pi i$.

\begin{remark}
Although we will not use it, the points of the form $(\x ,{\bf 0})\in \widehat{\cal B}-{\cal B}$ with $\x =(\l ,-\l ,\ldots ,\l ,-\l )$ for $\l >0$ small, represent boundary points of ${\cal W}$ corresponding to  Riemann minimal examples close to a stack of vertical catenoids.
\end{remark}

\begin{proposition}
\label{proposholomextensioncat}
Both $g$ and $\phi$ depend holomorphically on all parameters $(\x ,\y )$ in a neighborhood of $({\bf 0},{\bf 0})$(including at points of ${\cal B}\cup \widehat{\cal B}$).
\end{proposition}
\begin{proof}
Same as Proposition \ref{proposholomextensionscherk1p}.
\end{proof}

\begin{proposition}
\label{proposABextCat}
For $1\leq j\leq 2k$, the functions $A_j$ and $\widetilde{B}_j=y_jB_j$
extend holomorphically through ${\cal B}\cup \widehat{\cal B}$.
\end{proposition}
\begin{proof}
The extendability of $A_j$ through a point $(\x _0,\y _0)\in {\cal B}\cup \widehat{\cal B}$ is a consequence of  Proposition~\ref{proposholomextensioncat} since the curves $\Gamma_h$ stay in the limit Riemann
surface minus its nodes (if any).
For $\widetilde{B}_j$ we cannot apply Proposition~\ref{proposholomextensioncat} directly
because some of the points $O_{j-1},O_j$ or $\infty _{j-1},\infty _j$ where we compute the residues may collapse either into node points (when $(\x _0,\y _0)\in {\cal B}$) or into branch points of $g$ (when $(\x _0,\y _0)\in
\widehat{\cal B}-{\cal B}$). Rather we estimate the rate at which these residues blow-up.

We prove the extendability of $\widetilde{B}_j$ when $j$ is odd, the proof for $j$ even is similar.
Consider a tuple $(\x ,\y )\in D(0,\ve )^{4k}-({\cal B}\cup \widehat{\cal B})$ close to $({\bf 0},{\bf 0})$,
and let $\Omega \in g^{-1}(\overline{\C })$ be the annulus bounded by $\Gamma_{j-1}$ and
$\Gamma_j$, where $g=\aleph (\x ,\y )$. We shall first make a conformal representation of
this domain into a standard annulus.
$\Omega$ contains the branch points $a=a_j$, $b=b_j$. Recall
that $a$, $b$ are close to $0$.
We introduce the functions
\[
u=\sqrt{\frac{z-a}{z-b}},\qquad v=\frac{1+u}{1-u}
\]
and fix the sign of the square root by asking that $u\sim 1$ on
$\Gamma_j$ and consequently $u\sim -1$ on $\Gamma_{j-1}$. Then
 both $u$ and $v$ are conformal representations of $\Omega$.
On $\Gamma_j$ we may write $z=e^{i\theta}$
$$u=\sqrt{1+\frac{(b-a)e^{-i\theta}}{1-be^{-i\theta}}}
\sim 1+\frac{1}{2}(b-a)e^{-i\theta}$$
when $a,b$ are close to $0$.
Consequently, $u(\Omega)$ is close to $\overline{\C}$ minus the two disks
$D(\pm 1, r/2)$ where $r=|b-a|$, and $v(\Omega)$ is close to
the annulus $D(0,4/r)-D(0,r/4)$.
So we may write the Laurent series of $\phi$ with respect to the
variable $v$ in the annulus $D(0,1/r)-D(0,r)$ :
\[
\phi=\sum_{n\in\Z}c_nv^n dv
\]
where $c_n=\frac{1}{2\pi i}\int_{|v|=1}\frac{\phi}{v^{n+1}}$ depends holomorphically on all parameters.
Since $\phi$ extends holomorphically, the integral of $|\phi |$ on each unit circle $|z|=1$
extends analytically. Hence it is bounded by some constant $C>0$ (for the remainder of the proof
we shall denote by $C$ different positive constants). On $\Gamma_j$ we have $|v|\sim 4/r$,
hence $|v|\geq 1/r$. This gives if $n+1\geq 0$
\[
|c_n| =\frac{1}{2\pi}\left|\int_{\Gamma_j} \frac{\phi}{v^{n+1}}\right| \leq
\frac{1}{2\pi}\int_{\G_j}\frac{|\phi|}{|v|^{n+1}}\leq C r^{n+1}.
\]
On $\G_{j-1}$ we have $|v|\sim r/4$ hence $|v|\leq r$. This gives if $n+1\leq 0$
\[
|c_n|\leq \frac{1}{2\pi}\int_{\Gamma_{j-1}} |\phi| |v|^{|n+1|}\leq Cr^{|n+1|}.
\]
In summary,
\begin{equation}
\label{eq:estimccat}
|c_n|\leq C r^{|n+1|}\qquad \mbox{for all }n\in \Z .
\end{equation}
We now compute the residues appearing in $B_j$. First note that $z=0$ is equivalent to $v=\alpha$ or $v=1/\alpha$, where
\[
\alpha=\frac{1+\sqrt{a/b}}{1-\sqrt{a/b}}.
\]
Furthermore,  $g=z=\frac{a-b}{4v}(v-\alpha)(v-1/\alpha)$ hence
\[
\Res_{v=\alpha}\left( g^{-1}\phi\right) =\frac{4}{(a-b)(\alpha-1/\alpha)}\sum_{n\in\Z}
c_n\alpha^{n+1}
=-\frac{1}{\sqrt{ab}}\sum_{n\in\Z} c_n \alpha^{n+1}.
\]
(Observe here that this residue is indeed a multivalued function
as was explained in Remark~\ref{remarkResmultivalued}).
From the equality $(b-a)\alpha=(\sqrt{a}+\sqrt{b})^2$
we see that $|\alpha|\ll 1/r$ and $|\alpha|^{-1}\ll 1/r$. Using this information,
the estimate (\ref{eq:estimccat}) and the convergence of the series $\sum _{n\geq 1}t^n$ for any $t\in (0,1)$,
it is straightforward to check that $\sqrt{y_j} \Res (g^{-1} \phi)$ is bounded, where the residue is computed
at $0_j$ or $0_{j-1}$. Hence $y_j B_j$ is bounded on the set $y_j\neq x_j^2$, $y_j\neq 0$. Since
$y_j B_j$ is a well defined holomorphic function, it extends holomorphically by the Riemann extension theorem.
\end{proof}

\subsection{Partial derivatives and Inverse Function Theorem.}
\begin{proposition}
\label{proposAjBjderiv}
For each $j=1,\ldots ,2k$ it holds
\[
A_j({\bf 0},{\bf 0})=0,\qquad \widetilde{B}_j({\bf 0},{\bf 0})=-1,\qquad
\frac{\partial A_{j-1}}{\partial x_j}({\bf 0},{\bf 0})=2\pi i
\]
(with $A_0$ understood as $A_{2k}$), and all remaining partial derivatives of the $A_h$ are zero.
We will not need the partial derivatives of the $\widetilde{B}_h$.
\end{proposition}
\begin{proof}
Let $(\x ,\y)_j=(0,\ldots ,0,x_j,0,\ldots ,0,y_j,0,\ldots ,0)$ with $y_j\neq x_j^2$ and $y_j\neq 0$.
The Riemann surface associated to $(\x, \y )_j$ has $2k-1$ nodes which disconnect it into $2k-1$ genus zero components.
On $2k-2$ of these components (which we will call {\it simple spheres}), the corresponding meromorphic map $g_m$ is
the usual degree one $z$-map if we see the simple sphere as a copy of $\overline{\C }$, and the height differential
$\phi$ becomes $\frac{dz}{z}$. The remaining component $S$ is obtained from the copies $\overline{\C }_{j-1},
\overline{\C }_j$ (so we will call $S$ a {\it doubly sphere}).

First consider the case $j$ odd. Then $S$ can be parametrized by $\{ (z,w) \in \overline{\C }^2 \ | \
w^2=(z-a)(z-b)\} $ where $a+b=2x_j$, $ab=y_j$, and
the corresponding meromorphic map is $g(z,w)=z$. Hence, $w=\sqrt{(z-a)(z-b)}$ is well defined on $S$, where we convine
that the sign of the square root is determined by $w\sim z$ in $\overline{\C }_j$ and $w\sim -z$ in $\overline{\C }_{j-1}$. Since $\phi$ has simple poles at the nodes $\infty_{j-1}$ and $\infty_{j}$ with respective residues $1$ and $-1$,
and
\[
\Res_{\infty_j}\frac{dz}{w}=\Res_{z=\infty}\frac{dz}{z\sqrt{1-2x_j/z+y_j/z^2}}=-1,
\]
we conclude that $\phi=dz/w$. Thus,
\begin{equation}
\label{eq_Aj}
A_j(\x ,\y )_j=\int _{\G _j}\frac{dz}{zw}=-2\pi i \Res _{\infty _j}\frac{dz}{zw}=0
\end{equation}
and
\begin{equation}
\label{eq_Aj-1}
A_{j-1}(\x ,\y )_j=\int _{\G _{j-1}}\frac{z\, dz}{w}=2\pi i \Res _{\infty _{j-1}}\frac{z\, dz}{-z\sqrt{1-2x_j/z+y_j/z^2}}=
2\pi i x_j.
\end{equation}
These two equations imply that $A_j({\bf 0},{\bf 0})=0$ for every $j=1,\ldots ,2k$ (not necessarily odd).
Furthermore, given $h=1,\ldots ,2k$ it holds
$\frac{\partial A_j}{\partial x_h}({\bf 0},{\bf 0})=\lim _{y_h\to 0}\left.\frac{d}{dx_h}\right|_{x_h=0}A_j(\x, \y )_h$.
If $h\neq j-1,j$ then the integral in $A_j(\x ,\y )_h$ can be computed in a simple sphere which does not
depend on $x_h$, hence $\frac{\partial A_j}{\partial x_h}({\bf 0},{\bf 0})=0$. If $h=j$ then
(\ref{eq_Aj}) implies $\frac{\partial A_j}{\partial x_j}({\bf 0},{\bf 0})=0$. If $h=j-1$, then (\ref{eq_Aj-1})
gives $\frac{\partial A_{j-1}}{\partial x_j}({\bf 0},{\bf 0})=2\pi i$.

Concerning $\widetilde{B}_j$, we write
\[
\Res_{0_j}\left( g^{-1}\phi \right) =\frac{1}{\sqrt{y_j}},\qquad
\Res_{0_{j-1}}\left( g^{-1}\phi\right) =\frac{-1}{\sqrt{y_j}},
\]
hence $\widetilde{B}_j(\x ,\y )_j=-1$ for all $j$ odd.
The computations in the case $j$ even are similar, with the following modifications:
\[
w=\sqrt{\left(\frac{1}{z}-\frac{1}{a}\right)
\left(\frac{1}{z}-\frac{1}{b}\right)}
=\sqrt{\frac{1}{z^2}-\frac{2x_j}{z}+y_j},\qquad \phi=\frac{-dz}{z^2w},
\]
\[
A_j(\x ,\y )_j=\int_{\Gamma_j}\frac{-dz}{zw}=2\pi i\Res_{0_j}\frac{dz}{zw}=0,
\]
\[
A_{j-1}(\x ,\y )_j=\int_{\Gamma_{j-1}}\frac{-dz}{z^3w}=2\pi i\Res_{0_{j-1}}\frac{-dz}{-z^2\sqrt{1-2x_jz+y_jz^2}}
=2\pi i x_j,
\]
\[
\Res_{\infty _j}\left( g\phi \right) =\frac{1}{\sqrt{y_j}},\qquad
\Res_{\infty _{j-1}}\left( g\phi \right) =\frac{-1}{\sqrt{y_j}},
\]
from where one obtains similar conclusions to the case $j$ odd. The derivatives of the $A_h$
with respect to the variables $y_j$ clearly all vanish, and the proof is complete.
\end{proof}

Recall that given $j=1,\ldots ,2j$, we have $1/B_j=y_j/\widetilde{B}_j$. Since $\widetilde{B}_j$ extends holomorphically
to the polydisk $D(0,\ve )^{4k}$ (Proposition~\ref{proposABextCat}) and $\widetilde{B}_j({\bf 0},{\bf 0})=0$
(Proposition~\ref{proposAjBjderiv}), we deduce that $1/B_j$ extends holomorphically to $D(0,\ve )^{4k}$ for
$\ve >0$ small enough. Hence, the map $\Theta :D(0,\ve )^{4k}\to \C^{4k}$ given by
\begin{equation}
\label{eqTheta}
\Theta (\x ,\y )\left(A_1,\ldots ,A_{2k},1/B_1,\ldots ,1/B_{2k}\right)
\end{equation}
is holomorphic and $\Theta ({\bf 0},{\bf 0})=({\bf 0},{\bf 0})$.

\begin{theorem}
\label{thmuniqcat}
There exists $\ve >0$ small such that $\Theta $ restricts to $D(0,\ve )^{4k}$ as a biholomorphism
onto its image.
\end{theorem}
\begin{proof}
It follows from Proposition~\ref{proposABextCat} that
\[
\frac{\partial B_j^{-1}}{\partial x_i}({\bf 0},{\bf 0})=0,\qquad
\frac{\partial B_j^{-1}}{\partial y_i}({\bf 0},{\bf 0})=-\delta_{ij},
\]
Using these equations and the values of the partial derivatives of the functions $A_h$ at $({\bf 0},{\bf 0})$
given by Proposition~\ref{proposABextCat}, it is straightforward to check that the Jacobian matrix of
$\Theta $ at $({\bf 0},{\bf 0})$ is invertible. Now the Theorem follows from the Inverse function theorem.
\end{proof}

\begin{remark}
\label{remarklocaluniqcat}
As a consequence of Theorem~\ref{thmuniqcat}, for $t<0$ and $b\in \C$ close to $({\bf 0},{\bf 0})$
there exists a unique $(\x,\y)\in D(0,\ve )^{4k}$ such that $\Theta (\x,\y)=(b,\overline{b},\cdots,b,\overline{b},
t,\cdots,t)$. This uniqueness result and equation (\ref{eq:gclosesperiods}) imply that the space of immersed
minimal surfaces sufficiently close to the boundary point of ${\cal W}$ given by the $2k$ catenoids
(see Proposition~\ref{propos2kcatenoids}) has three real freedom parameters.
\end{remark}

\section{Openness.}
\label{secopen}

Recall that ${\cal K}\subset {\cal S}$ represents the space of
standard examples with $4k$ ends. A direct consequence of its construction is that ${\cal K}$ is
closed in ${\cal S}$. We saw in Section~\ref{secstandardexamples} that ${\cal K}$ is open in ${\cal S}$, by
the nondegeneracy of any standard example. Both closeness and openness remain valid for the space $\widetilde{\cal K}$
of marked standard examples inside $\widetilde{\cal S}$. By Lemma~\ref{lemmaJ}, $\widetilde{\cal K}$ is closed and open in ${\cal M}$. Our Theorem~\ref{thm1} reduces to prove that $\widetilde{\cal K}=\widetilde{\cal S}$,
which will be proved by contradiction in Section~\ref{secproofmainthm}. For this reason, in what follows we will assume $\widetilde{\cal S}-\widetilde{\cal K}\neq \mbox{\O }$.

\begin{theorem}\label{thopen}
The classifying map $C:\widetilde{\cal S}-\widetilde{\cal K}\to\R^*\times \C$ is open.
\end{theorem}
\begin{proof}
Fix a marked surface $M\in\widetilde{\cal S}-\widetilde{\cal K}$. It suffices to see that $C$ is open in a neighborhood of $M$ in $\widetilde{\cal S}-\widetilde{\cal K}$. Let $(a,b)=C(M)\in\R^*\times \C$ and ${\cal M}(a,b)=L^{-1}(L_{(a,b)})$ $\subset {\cal M}$ (with the notation of Subsection~\ref{subsecligaturemap}).
Since $\widetilde{\cal K}$ is open and closed in $\widetilde{\cal S}$ and $\widetilde{\cal S}(a,b)=\widetilde{\cal S}\cap{\cal M}(a,b)$ is an analytic subvariety of ${\cal W}$ by Proposition~\ref{proposanalsubv}, we conclude that
${\cal (\widetilde{S}-\widetilde{K})}(a,b)={\cal (\widetilde{S}-\widetilde{K})}\cap \widetilde{\cal S}(a,b)$ is an analytic subvariety of ${\cal W}$.

\begin{assertion}
\label{asSKcomp}
${\cal (\widetilde{S}-\widetilde{K})}(a,b)$ is compact.
\end{assertion}
To prove Assertion~\ref{asSKcomp}, take a sequence
$\{M_n\}_n\subset {\cal (\widetilde{S}-\widetilde{K})}(a,b)$. By Proposition~\ref{proposcurvestim}, the
sequence of geometric surfaces $\{ M_n\} _n$ has uniformly bounded Gaussian curvature. Similarly as in the
proof of Proposition~\ref{propos2kScherks1p}, we can find $\de >0$ such that for every $n$, there exists
a point $p_n\in M_n$ where the normal vector $N_n(p_n)$ lies in $\esf ^2\cap \{ x_2=0\} $ and stays at spherical distance at least $\de $ from the branch locus of the Gauss map $N_n$ of $M_n$. The usual limit process shows that, after passing to a subsequence, suitable liftings of $M_n-p_n$ converge smoothly as $n\to \infty$ to a properly embedded nonflat minimal surface $\widetilde M_\infty\subset \R^3$ in one of the six cases listed in
Proposition~\ref{proposdifferentlimits}.
As $\widetilde{\cal K}$ is open in $\widetilde{\cal S}$, if $\widetilde{M}_{\infty }$ were in case
{\it (vi)} of Proposition~\ref{proposdifferentlimits} then its quotient would actually be in
$\widetilde{S}-\widetilde{K}$, hence it only remains to discard the first five possibilities of Proposition~\ref{proposdifferentlimits} for $\widetilde M_\infty$.

As the period vector $H_n$ at the ends of $M_n$ is $(0,\pi a,0)$ for all $n$, the cases
{\it (i),(ii),(iii)} of Proposition~\ref{proposdifferentlimits} are not possible for
$\widetilde M_\infty$. If $\widetilde M_\infty$ is a singly periodic
Scherk surface $S_\rho$ with two horizontal ends, then any choice of the nonhorizontal period
vector $T_n$ of $M_n$ necessarily diverges to $\infty$ (case {\it (iv)} of
Proposition~\ref{proposdifferentlimits}).
By Proposition~\ref{propos2kScherks1p} and Remark~\ref{notalimscherk1p}-{\it (ii)},
for all $n$ large we can find a new marked surface $M'_n$ which can be
viewed inside the open subset ${\cal U}(\ve')$  of $\C^{4k}$ appearing in
Theorem~\ref{teorLbiholScherk1p}, such that
$M_n,M'_n$ only differ in the homology class of the last component of the
marked surfaces. Since $L|_{{\cal U}(\ve')}$ is a biholomorphism (Theorem~\ref{teorLbiholScherk1p}),
the space of tuples in ${\cal U}(\ve ')$ producing immersed minimal surfaces has three real freedom
parameters. Since $\widetilde{\cal K}$ has real dimension three and $C|_{\widetilde{\cal K}}$ takes
values arbitrarily close to $(a,b)$, we deduce that
if $M'\in \widetilde{\cal S}$ lies in ${\cal U}(\ve')$, then $C(M')$ coincides
with the value of $C$ at a certain standard example $M_0\in\widetilde{\cal K}$. In particular
$L(M')=L(M_0)$ and as $L|_{{\cal U}(\ve')}$ is injective, we have $M'=M_0$. As this last equality cannot hold for $M'$ being $M_n'\in
\widetilde{\cal S}-\widetilde{\cal K}$, we conclude that $\widetilde M_\infty$ is not a singly periodic
Scherk minimal surface.

Now assume that $\widetilde M_\infty$ lies in case {\it (v)} of Proposition~\ref{proposdifferentlimits}.
Let $\G$ be a component of the intersection of $\widetilde M_\infty$
with a horizontal plane $\{ x_3=c\}$ whose height does not coincide with the
heights of the horizontal ends of $\widetilde M_\infty$. Since $\widetilde{M}_{\infty }$
has exactly two nonhorizontal ends, $\G$ is an embedded
$U$-shaped curve with two almost parallel divergent ends, and
if we denote by $\Pi \subset \R^3$ the plane passing
through the origin parallel to the nonhorizontal ends of $\widetilde M_\infty$, then
the conormal vector to $\widetilde M_\infty$ along each of the divergent branch of $\G$ becomes
arbitrarily close to the upward pointing unit vector $\eta \in \Pi$ such that $\eta
$ is orthogonal to $\Pi \cap \{ x_3=c\} $. Since translated liftings of the $M_n$
converge smoothly to $\widetilde M_\infty$, we deduce that $M_n$ contains
arbitrarily large arcs at constant height along which the conormal vector $\eta_n$ is
arbitrarily close to $\eta $. In particular, the integral of the third component of
$\eta_n$ along such arcs becomes arbitrarily large. As the conormal vector of
$M_n$ along any compact horizontal section misses the horizontal values (the Gauss map of $M_n$ is never vertical),
it follows that the vertical component
of the flux of $M_n$ along a compact horizontal section diverges to $\infty$ as $n\to
\infty $. This contradicts our normalization on the surfaces of ${\cal S}$, and
proves Assertion~\ref{asSKcomp}.
\par
\vspace{.2cm}

We now finish the proof of Theorem~\ref{thopen}. By Assertion~\ref{asSKcomp} and Lemma~\ref{lemacptanalsub}, ${\cal (\widetilde{S}-\widetilde{K})}(a,b)$ is a finite subset, hence we can
find an open set ${\cal U}$ of ${\cal W}$ containing $M$ such that  ${\cal (\widetilde{S}-\widetilde{K})}(a,b)\cap
{\cal U}= {\cal
M}(a,b)\cap {\cal U}=\{ M\} $. In terms of the ligature map $L:{\cal W}\to \C^{4k}$,
the last equality writes as $L^{-1}(L_{(a,b)})\cap {\cal U}=
\{ M\} $. Since $L$ is holomorphic, we can apply the Openness Theorem for finite holomorphic
maps (see~\cite{GriHar} page 667) to conclude that $L|_{\cal U}$ is an open map.
Finally, the relationship between the ligature map $L$ and the map $C$ gives
the existence of a neighborhood of $M$ in ${\cal \widetilde{S}-\widetilde{K}}$ where the restriction of
$C$ is open.
\end{proof}

The same argument in the proof of Assertion~\ref{asSKcomp} remains valid under the
weaker hypothesis on $C(M_n)$ to converge to some $(a,b)\in\R^*\times \C$ instead of being constant on a sequence
$\{M_n\}_n\subset{\cal \widetilde{S}-\widetilde{K}}$. This proves the validity of the following statement.

\begin{theorem}\label{thproper}
The classifying map $C:{\cal \widetilde{S}-\widetilde{K}}\to\R^*\times \C$ is proper.
\end{theorem}

\section{The proof of Theorem~\ref{thm1}.}
\label{secproofmainthm}
Recall that we were assuming ${\cal \widetilde{S}-\widetilde{K}}\neq \mbox{\O }$. By Theorems~\ref{thopen} and \ref{thproper},
$C:{\cal \widetilde{S}-\widetilde{K}}\to\R^*\times \C$ is an open
and proper map. Thus, $C({\cal \widetilde{S}-\widetilde{K}})$ is an open and closed subset of
$\R^*\times \C$. Since $C({\cal \widetilde{S}-\widetilde{K}})$ has points in both connected components of
$\R^*\times \C $, we deduce that $C|_{\widetilde{\cal S}-\widetilde{\cal K}}$ is surjective. In particular,
we can find a sequence $\{M_n\}_n\subset {\cal \widetilde{S}-\widetilde{K}}$ such that $\{C(M_n)\}_n$ tends to
$(\infty,0)$ as $n$ goes to infinity. Now the argument is similar to the one in the proof of Assertion~\ref{asSKcomp}
when we discarded the singly periodic Scherk limit, using Proposition~\ref{propos2kcatenoids}, Remark~\ref{notalimcat}
and Theorem~\ref{thmuniqcat} instead of Proposition~\ref{propos2kScherks1p}, Remark~\ref{notalimscherk1p}-{\it (ii)} and
Theorem~\ref{teorLbiholScherk1p}.

\center{Joaqu\'\i n P\' erez at jperez@ugr.es\\
Magdalena Rodr\'\i guez at magdarp@ugr.es\\
Martin Traizet at traizet@univ-tours.fr}

\end{document}